\documentclass[12pt]{article}

\usepackage{amssymb}
\usepackage{amsmath}
\usepackage{enumerate}
\usepackage{graphicx}
\usepackage{amscd}
\usepackage{latexsym}
\usepackage{epsfig}
\usepackage[latin1]{inputenc}
\usepackage{hyperref}

\usepackage{a4wide}

\def\h{\mathbb{H}}
\def\r{\mathbb{R}}
\def\n{\mathbb{N}}

\def\c{\mathbb{C}}
\def\s{\mathbb{S}}
\def\d{\mathbb{D}}
\def\l{\mathbb{L}}
\def\z{\mathbb{Z}}

\def\rb{\mathcal{R}}
\def\sb{\mathcal{S}}

\def\mb{\mathcal{M}}

\def\pg{\mathfrak{p}}

\def\Gg{\mathfrak{G}}

\def\pg{\mathfrak{p}}

\def\sg{\mathfrak{s}}

\def\mg{\mathfrak{m}}
\def\Mg{\mathfrak{M}}

\def\vg{\mathfrak{v}}

\newcommand{\df}{ \stackrel{\rm def}{=}}

\def\div{\mathfrak{Div}}
\def\coc{\l^3/_{\langle (1,0,0)\rangle}}
\def\nb{\mathcal{N}}
\def\N{\mathcal{N}}

\def\T{\cal T}
\def\E{\cal E}
\def\J{\cal J}

\def\V{{\cal V}}

\newtheorem{lemma}{Lemma}[section]
\newtheorem{remark}{Remark}[section]
\newtheorem{theorem}{Theorem}[section]
\newtheorem{proposition}{Proposition}[section]
\newtheorem{corollary}{Corollary}[section]

\newtheorem{definition}{Definition}[section]
\newenvironment{proof}{\trivlist
\item[\hskip\labelsep{\em Proof}\,:]}{\hfill{$\Box$}\endtrivlist}

\title{\LARGE The moduli space of embedded singly periodic  maximal surfaces with isolated
 singularities in the Lorentz-Minkowski space $\l^3$}
\author{\Large Isabel Fern\'{a}ndez   \thanks{Research partially
supported by
MCYT-FEDER grant number MTM2004-00160. \newline 2000 Mathematics
Subject Classification. Primary 53C50; Secondary 58D10, 53C42.
 }
  \and \Large  Francisco J. L\'{o}pez $ ^{{\ast} }$
  \and \Large  Rabah Souam  $ ^{{\ast} }$}

\begin{document}

\maketitle

\begin{abstract}
 We show that, up to some natural normalizations, the moduli space
 of singly periodic complete embedded maximal surfaces in the
 Lorentz-Minkowski space $\l^3=(\r^3,dx_1^2+dx_2^2-dx_3^2),$ with fundamental piece having a finite number
 $(n+1)$ of singularities,  is a real  analytic manifold of
 dimension $3n+4.$ The underlying topology agrees with the topology of
 uniform convergence of graphs
 on compact subsets of $\{x_3=0\}.$
 \end{abstract}

\section{Introduction} \label{sec:intro}

\begin{figure}
  \centering
  \includegraphics[width=5cm]{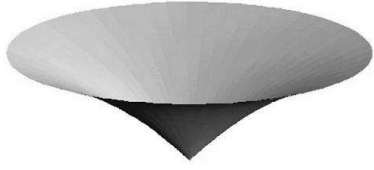}
  \caption{The Lorentzian catenoid}\label{fig:cat}
\end{figure}

A maximal surface  in the Lorentz-Minkowski space $\l^3=(\r^3,
dx_1^2+dx_2^2-dx_3^2)$ is a spacelike surface with zero mean curvature. It locally
maximizes   the area functional associated to variations by spacelike surfaces.  In
a  pioneering work, Calabi \cite{calabi}  proved  that the  affine spacelike planes
are the only  complete maximal surfaces in   $\l^3$ (Calabi, in fact,  showed the
analogous result for maximal hypersurfaces in $\l^4$ and this was later extended to
maximal hypersurfaces in $\l^n,$ for all $n,$ by Cheng and Yau \cite{cheng-yau}).
Nonetheless,  recent works show there is a rich global theory of complete maximal
surfaces with singularities in $\l^3$ (\cite{fer-lop-sou},  \cite{fer-lop},
\cite{klyachin}, \cite{lopez-lopez-souam}, \cite{um-ya}). For instance, Umehara and
Yamada \cite{um-ya} have obtained results on the global behavior of immersed maximal
surfaces having analytic curves of singularities. Of particular interest are the
complete embedded maximal surfaces having a closed discrete set of singularities.
First it should be remarked that maximal surfaces in   $\l^3$ share some properties
with minimal surfaces in the Euclidean space $\r^3.$ Indeed,  Kobayashi
\cite{kobayashi} gave  a Weierstrass type representation in terms of meromorphic
data similar to the one of minimal surfaces in $\r^3.$ Also, both types of surfaces
are  locally represented as graphs of solutions of elliptic operators.     An
important difference, however,  is that the maximal surface equation for graphs in
$\l^3$  may have solutions with isolated singularities and this never happens for
the minimal surface equation in $\r^3.$  Otherwise said, an embedded maximal surface
can  have an isolated singularity, contrarily to an embedded minimal surface in
$\r^3.$  This is illustrated by the Lorentzian half-catenoid, \cite{kobayashi} (see Figure \ref{fig:cat}).  At
a singular point around which  a maximal surface  is embedded  in $\l^3,$ the Gauss
curvature blows up, the limit tangent planes become lightlike and the surface is
asymptotic to a half lightcone at the singularity (cf. \cite{klyachin},
\cite{kobayashi} and \cite{fer-lop-sou}). For these reasons, such  points are called
{\it conelike} singularities.

Embedded complete maximal surfaces with  a closed discrete set of
singularities  are global graphs over any spacelike plane.
Conformally, the regular set of such a  surface  is a Riemann
surface minus as many pairwise disjoint closed conformal disks
(without accumulation points) as singular points in the surface. A
fundamental  observation is that the Weierstrass data for the
surface extends to the double Riemann surface, \cite{fer-lop}.
This allows one to work on boundaryless Riemann surfaces.

In \cite{fer-lop-sou} we developed the theory of embedded complete
maximal surfaces in $\l^3$ of finite type, that is those having a
finite number of singularities. We showed, in particular, that the
moduli space of such surfaces with $n+1\ge 2$ singularities and
vertical limit normal at infinity is a $3n+4-$dimensional
manifold.

The next simplest subclass of maximal surfaces in $\l^3$ with a
closed discrete set of singularities consists of the {\it
periodic} ones with finite type in the quotient. This refers to
surfaces that are invariant under a discrete group $G$ of isometries of
$\l^3$ acting freely and properly and such that the quotient surfaces are
embedded and have a finite number of singularities in the flat and
complete 3-dimensional Lorentzian manifold $\l^3/G.$ In
\cite{fer-lop}, the first and second author proved fundamental
facts about the global geometry of these surfaces. In particular,
they classified the discrete groups $G\subset Iso(\l^3)$ for which $\l^3/G$
contains complete maximal surfaces of finite type. If $\l^3/G$
is orientable and orthochronous (i.e the elements of $G$ preserve
the orientation and the future time direction) and contain a
complete embedded maximal surface $S$ of finite type then $G$ is a group of
spacelike translations of rank one or two  and $S$ is an annulus
of finite conformal type or a torus, respectively.

\begin{figure}
  \centering
  \includegraphics[width=11cm]{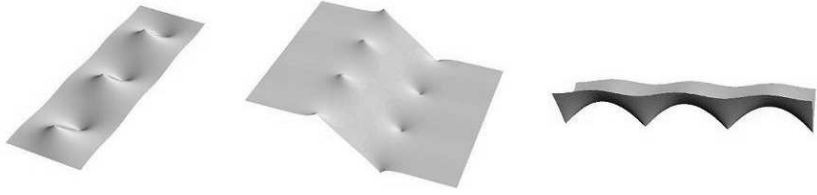}
  \caption{Examples of CMSF surfaces}\label{fig:scherk}
\end{figure}

In this paper, we study the moduli space of (embedded) {\it singly periodic}
maximal surfaces in $\l^3$ having  finite type in the quotient
space. That is, we consider the group $<T>$ generated by a
spacelike translation $T$ in $\l^3$  and complete embedded maximal surfaces
with a finite number of singularities in $\l^3/<T>.$ Several examples
of this kind were constructed in
\cite{fer-lop} (see figure \ref{fig:scherk}).

First note that up to an ambient
isometry of $\l^3$ and rescaling, we can suppose that $T=(1,0,0).$
In this case any complete embedded maximal surface with a finite number of singularities in $\coc$
is a graph over the cylinder $\{x_3=0\}\subset\coc.$ We also normalize so that one of the ends of the surface is asymtpotic to $\{x_3=0,x_2\ge 0\}.$
Our main result then says that:

\begin{quote}{\em  The space $\Mg_n$ of marked entire maximal
graphs over the cylinder $\{x_3=0\}$ in $\l^3/<(1,0,0)>$ having
$n+1 \geq 2$  conelike singularities  (the mark is an ordering of
the set of singularities)  and an end asymptotic to $\{x_3=0,
x_2\ge 0\},$ is a real analytic manifold of dimension $3n+4. $
A global coordinate system is given by the ordered sequence of
points in the mark and the normal to the second end. This space is
a $(n+1)!$-sheeted covering of the space $\Gg_n$ of (non marked)
entire maximal graphs over the cylinder $\{x_3=0\}$ in
$\coc,$  having  $n+1 \geq 2$  conelike singularities
and  an end asymptotic to $\{x_3=0, x_2\ge 0\}$. The underlying
topology of $\Gg_n$ is equivalent to the uniform convergence of
graphs on compact subsets  of the cylinder $\{ x_3=0\}.$}
\end{quote}

We have organized our paper as follows:  Section 2 contains  some preliminaries
about the local behavior of maximal surfaces around isolated singularities and the
global behavior of complete maximal surfaces of finite type in the quotient space
$\l^3/<T>.$ Section 3 is devoted to the proof of the main theorem. Our approach relies on
algebraic geometry tools:
we define some natural bundles on the moduli space ${\cal
T}_n$ of  once punctured marked circular domains with $n+1$ boundary components (a
mark  is an ordering of the boundary circles), and introduce a spinorial bundle
${\cal S}_n$ associated to the moduli space of Weierstrass data of surfaces in the
space of graphs with $n+1$ singularities. The convergence in $\Mg_n$ means convergence of
marked  conformal structures in ${\cal T}_n$ and of Weierstrass data.

\section{Preliminaries} \label{sec:prelim}

We denote by $\overline{\c},$ $\d$ the extended complex plane
$\c \cup \{\infty\}$ and the unit disc $\{z \in \c \;:\; |z|<1\},$ respectively.

Throughout this paper, $\l^3$ will denote the three dimensional
Lorentz-Minkowski space $(\r^3,\langle , \rangle),$ where $\langle ,
\rangle=dx_1^2+dx_2^2-dx_3^3.$ By definition, a coordinate system $(y_1,y_2,y_3)$ in $\l^3$  is said to be a {\em $(2,1)$-coordinate system} if the Lorentzian metric is given by $dy_1^2+dy_2^2-dy_3^3.$
We say that a vector ${\bf u} \in \r^3- \{ {\bf 0} \}$ is spacelike,
timelike or lightlike if
$\|u\|^2:=\langle {\bf u}, {\bf u} \rangle$ is positive, negative or zero,
respectively. When $u$ is spacelike, $\|u\|$ is chosen non negative. The vector
${\bf 0}$ is spacelike by definition.  A plane in $\l^3$ is spacelike,
timelike or lightlike if the induced metric is Riemannian, non degenerate
and indefinite or degenerate, respectively.

We call $\h^2 = \{ (x_1,x_2,x_3) \in \r^3 \;:\;
x_1^2+x_2^2-x_3^2=-1\}$ the hyperbolic sphere in $\l^3$ of constant
intrinsic curvature $-1.$ Note that $\h^2$ has two connected  components $\h^2_+:=\h^2 \cap \{x_3 \geq 1\}$ and $\h^2_-:=\h^2 \cap \{x_3 \leq -1\}.$ The stereographic
projection $\sigma$ for $\h^2$ is defined as follows:
$$\sigma:\overline{\c} - \{|z|=1\} \longrightarrow \h^2 \,; \; z \rightarrow
\left(\frac{2 \mbox{Im} (z)}{|z|^2-1}, \frac{-2 \mbox{Re} (z)}{|z|^2-1},
\frac{|z|^2+1}{|z|^2-1} \right),$$ where $\sigma(\infty)=(0,0,1).$

By definition, an isometry in $\l^3$ is said to be orthochronous if its associated linear isometry  preserves $\h^2_+$ (and so $\h^2 _-$). In other words, it preserves the future direction.\\

In the sequel, $\N$ will denote a complete flat 3-dimensional Lorentzian manifold (i.e., a 3-dimensional differential manifold endowed with a flat metric of index one). It is well known that the universal isometric covering of $\N$ is $\l^3$ (see for example \cite{oneill},\cite{wolf}).Thus  $\N$ can be regarded as the quotient of $\l^3$ under the action of a discrete group $G$ of isometries acting freely and properly on $\l^3.$

In what follows,  $\mb$ will denote a differentiable surface.

An immersion $X:\mb \longrightarrow \N$ is spacelike if
the tangent plane at any point is spacelike, that is to say, the induced metric on $\mb$ is Riemannian.
In this case, $S=X(\mb)$ is said to be a spacelike surface in $\N.$ If $\N=\l^3/G,$ where $G$ is a (possibly trivial) group of translations acting freely and properly on $\l^3,$ the locally well defined Gauss map $N_0$ of $X$ assigns to each point of $\mb$ a point of $\h^2.$  A connectedness argument gives that $N_0$ is globally well defined and $N_0(\mb)$ lies, up to a Lorentzian isometry, in $\h^2_-.$ This means that $\mb$ is orientable.

A maximal immersion $X:\mb \longrightarrow \N$ is a spacelike immersion
with null mean curvature. In this case, $S=X(\mb)$ is said to be a maximal surface in $\N.$
Using isothermal parameters, $\mb$ can be endowed with a conformal structure. In the orientable case, $\mb$ becomes a Riemann surface.


\begin{theorem}[Weierstrass representation of maximal surfaces in $\l^3$ \cite{kobayashi}]

Let $X:\mb\to\l^3$ be a conformal maximal immersion of a Riemann
surface. Then $g\df \sigma^{-1}\circ N_0$ is a meromorphic function, and there exists a
holomorphic 1-form $\phi_3$ defined on $\mb$ such that
\begin{enumerate}[(i)]
\item the 1-forms given by $\phi_1=\frac{i}{2}(\frac{1}{g}-g)\phi_3$
and $\phi_2=\frac{-1}{2}(\frac{1}{g}+g)\phi_3$ are holomorphic on
$\mb$

\item the induced Riemannian metric on $M$ is given by
$ds^2=|\phi_1|^2+|\phi_2|^2-|\phi_3|^2=\frac{1}{4}\big(\frac{1}{|g|}-|g|)^2|\phi_3|^2$

\item for any closed curve $\gamma$ in $M$ we have
$\mbox{Re}\int_\gamma(\phi_1,\phi_2,\phi_3)=0$

\item up to a translation, the immersion is given by
$X=\mbox{Re}\int_{P_0} (\phi_1,\phi_2,\phi_3),$ where $P_0\in M$
is an arbitrary point.
\end{enumerate}

Conversely, given $g$ and $\phi_3$ a meromorphic function and a
holomorphic 1-form on $M,$ respectively, such that $(i),$ $(ii)$ and
$(iii)$ are satisfied, then $(iv)$ defines a conformal  maximal
immersion of $M$ in $\l^3.$
\end{theorem}

\begin{remark}[Weierstrass data of maximal surfaces in translational spaces]
The Weierstrass data $(\phi_1,\phi_2,\phi_3)$ of a maximal surface
in $\l^3$ are invariant by translations. Therefore, maximal
surfaces in a quotient $\l^3/G,$ where $G$ is a group of
translations acting properly and freely, also have a Weierstrass
representation as above except that the condition $(iii)$ is replaced by the following one:

$(iii')$ for any closed curve $\gamma$ in $M$  the translation of associated
   vector
$\mbox{Re}\int_\gamma(\phi_1,\phi_2,\phi_3)$ is an element of the group $G.$
\end{remark}

\begin{definition}
Let $X: \mb \to \nb$ be a topological embedding of a smooth
surface in a Lorentzian 3-manifold and $F\subset \mb$ a closed
discrete subset. We say that $X$ is a maximal embedding
with singular set $F$ if $X|_{\mb- F}$ is a maximal embedding and
the induced metric on $\mb-F$ converges to zero at any point of $F.$

In this case we say that $S=X(\mb)$ is a maximal surface in $\nb$ with singularities at $X(F).$
\end{definition}

\begin{lemma}[structure of embedded singularities in $\l^3$ \cite{fer-lop-sou}, \cite{ecker}] \label{lem:sing}

Let  $X: {\cal D}\to \l^3$ be a maximal embedding defined on an
open  disk ${\cal D},$ with an isolated singularity at $q \in
{\cal D}.$

Then ${\cal D}- \{q\}$ is conformally equivalent to $\{z\in
\c\, :\, r<|z|< 1\},$ for some $r>0,$ and if $X_0:\{ r<|z|< 1\} \to \l^3 $ is
a conformal reparameterization of $X,$ then $X_0$ extends
analytically to $A_r:= \{1\leq |z|<1\}$ by setting $X_0(\{|z|=1\})=P_0:=X(q).$

The Weierstrass data $(g,\phi_3)$ of $X_0$ satisfy: $g$ is injective and $|g|=1$ on $\{|z|=1\},$  and $\phi_3(z)\ne 0,\,
|z|=1.$
In particular, $X_0$  reflects analytically about $\{|z|=1\}$ to the mirror surface $A_r^*:=\{z\in\c\; :\; 1 \leq |z|< 1/r\},$ verifying
$g\circ J=1/\overline{g}$ and $J^*(\phi_3)=-\overline{\phi}_3,$ where $J(z)=1/\overline{z}$ is the mirror involution.

Moreover for any spacelike plane  $\Pi$  plane containing $P_0$
the Lorentzian orthogonal projection  $\pi: X({\cal D}) \to
\Pi$ is a local homeomorphism and $X({\cal D})$ is asymptotic near $P_0$
to a half
light cone with vertex at $P_0.$

The point $P_0$ is said to be a conelike singularity of $X({\cal D}).$

\end{lemma}

\begin{remark}

The universal covering of a complete flat Lorentzian 3-dimensional  manifold
is isometric to $\l^3$ (cf. \cite{wolf}). Therefore the above
Lemma extends to the more general context of
complete flat Lorentzian 3-manifolds.
\end{remark}

As a consequence of the previous lemma,
if $X:\mb \to\nb$ is a complete embedded
maximal surface with a closed discrete set
$F\subset \mb$ of singularities, where $\nb$ is complete and flat,  then $\mb-F$ is conformally
equivalent to $\Sigma-\cup_{p\in F} D_p,$ where $\Sigma$ is a
Riemann surface and the $\{D_p\}_{p\in F}$ are closed pairwise
disjoint conformal disks with no accumulation in $\Sigma.$

We also have that the conformal reparameterization
$X:\Sigma-\cup_{p\in F}D_p\to \nb$ extends
analytically to $ \mb_0:=\Sigma-\cup_{p\in F}\mbox{Int}(D_p),$ by
putting $X(\partial D_p)=X(p)$ for each $p\in F.$ In the sequel we
will refer to $\mb_0$ as the {\em conformal support} of the
embedding $X.$
We also say that $\mb_0$ is the conformal support of the maximal surface $X(\mb)\subset\nb.$

In particular, if $\nb=\l^3$ or $\nb=\l^3/G,$ where $G$ is a
translational group, the Weierstrass data $(\phi_1,\phi_2,\phi_3)$
extend analytically to $\mb_0.$

We denote by $\mb_0^*$ the mirror surface of $\mb_0$ and by $N$ the double
surface, that is, $N=\mb_0\cup \mb_0^*$ with the identification
$\partial(\mb_0)\equiv \partial(\mb_0^*).$ Moreover, we call
$J:N\to N$ the antiholomorphic mirror involution, and observe that the
fixed point set of $J$ coincides with $\partial(\mb_0).$
By Lemma \ref{lem:sing}, the Weierstrass data $\Phi:=(\phi_1,\phi_2,\phi_3)$ can be extended by Schwartz reflection to $N$ satisfying
$J^*(\Phi)=-\overline{\Phi}.$\\

For the sake of simplicity, complete maximal embedded surfaces with a finite
set of singularities in $\nb$ will be called {\bf CMF} surfaces.

\begin{definition}[flux of a closed curve]\label{def:flux}
Let $S\subset\nb$ be an oriented CMF
surface, and let $X:\mb_0 \to \nb$ be a conformal
reparameterization of $S.$

For any closed curve $\gamma(s)$ in $\mb_0$ parameterized by the
arclength, we label $\nu$  as its unit conormal vector so that
$\{\nu,\gamma'\}$ is positive with respect to the orientation in
$\mb_0.$ The flux vector of the curve $\gamma$ is defined as
$$ F(\gamma):= \int_\gamma \nu(s)ds $$
\end{definition}

Since $X$ is harmonic it follows from Stokes theorem
that $F(\gamma)$ depends only on the homology class of $\gamma$ in
$\mb_0.$ If $\nb=\l^3$ or $\l^3/G,$ where $G$ is translational, and we denote
by $\Phi=(\phi_1,\phi_2,\phi_3)$ the Weierstrass data of $X,$ it is easy to check that
$$F(\gamma)= \mbox{Im} \big( \int_\gamma \Phi \big)$$

We define the flux at a conelike singularity $q\in S$ as the flux
along any curve homotopic to the boundary component $\gamma_0$ of
$\mb_0$ corresponding to $q.$ It can be checked hat the flux at conelike singularity is always a timelike
vector (see \cite{tubes}).

\begin{definition}[CMSF surfaces]
A complete embedded maximal surface $\tilde S\subset\l^3$ with a closed discrete set of
singularities is said to be singly periodic of finite
type if:
\begin{itemize}
\item $\tilde S$  is invariant under the
free and proper action  of an infinite cyclic group $G$ of isometries of $\l^3$
\item   $\tilde{S}/G$ is a CMF surface in $\l^3/G.$
\end{itemize}

In the sequel $\tilde{S}$ will be called for short a {\bf CMSF}
surface.

\end{definition}

\begin{remark} \label{re:conexion}
There is a natural connection between CMSF surfaces in $\l^3$ and CMF surfaces in  quotients $\l^3/G,$ where $G$ is cyclic.

As stated in the above definition a CMSF surface determines a CMF surface in the corresponding space $\l^3/G.$
Conversely, Theorem \ref{th:representation} below will show that the universal covering of a CMF surface in $\l^3/G$ is a CMSF surface in $\l^3.$
\end{remark}

\subsection{Representation of CMSF surfaces}

In the sequel we will denote by $\langle u\rangle$ the cyclic
group generated by the translation of vector $u\in\l^3.$

\begin{theorem}[\cite{fer-lop}]\label{th:representation}

Let $S$ be a CMF surface with $n+1$ singular
points in $\l^3/G,$ where $G$ is a cyclic group of isometries of $\l^3$ acting properly and freely. Then

\begin{itemize}

\item The CMSF surface $\tilde S$ obtained by lifting $S$ to $\l^3$ is an entire graph over any spacelike plane.

\item The group $G$ is generated by a spacelike translation $T.$

\item The conformal support $\mb_0$ of $S$ is $\c^{*}-\cup_{i=0}^{n}\mbox{Int}(D_i)$
where the $D_i,\, i=0,\ldots,n$ are pairwise disjoint closed
(circular) disks. The associated double surface $N$ is
$\overline{N}-\{0,\infty,J(0),J(\infty)\},$ where $\overline{N}$
is a compact Riemann surface of genus $n,$ the points $0,\infty\in
\overline{N}$ correspond to the ends of $\mb_0$ and $J$ denotes the
mirror involution.
\item If $X:\mb_0 \to \l^3/G$ is a conformal parameterization of $S,$ $X$ applies
each boundary circle $\partial(D_i)$ to a singular point of $S.$
\item The Weierstrass data $\Phi:=(\phi_1,\phi_2,\phi_3)$ of $X$  can be extended  by Schwarz reflection to $N,$
satisfying: $J^{*}(\Phi) = -\overline{\Phi}.$

Moreover $\Phi$ has simple poles at $0,\infty,J(0)$ and $J(\infty),$ and the topological ends of
$S$ are of Scherk type, that is to say, asymptotic to spacelike
flat half cylinders in $\l^3/G.$

\item $T$ can be chosen as the translation of vector $u=\mbox{Re}\int_{\gamma}\Phi,$
where $\gamma\subset \mb_0$ is a closed loop around $0.$

\end{itemize}

Conversely given $\mb_0:=\c^{*}-\cup_{i=0}^{n}\mbox{Int}(D_i)$ where
the $D_i,\, i=0,\ldots,n$ are pairwise disjoint closed (circular)
disks, define $N$ and $ J$ as before and take a Weierstrass data
$\Phi$ on $N$ satisfying $J^{*}(\Phi) = -\overline{\Phi}$ and
having simple poles at the ends. Then
$$X: \mb_0 \to \l^3/_{\langle u\rangle} \qquad
X=\mbox{Re}\int_{P_0}\Phi$$
where $u=\mbox{Re}\int_{\gamma}\Phi,$ defines a complete embedded
maximal surface with $n+1$ singularities and its universal
covering is a CMSF surface in $\l^3$ invariant under the group
$\langle u\rangle.$

\end{theorem}

\subsection{Uniqueness of  CMSF surfaces}

\begin{theorem}[Uniqueness]\label{th:uniqueness}

Let $S_i$ denote a CMF surface in
$\l^3/_{\langle u\rangle}$ with singular points $q_1^i,\ldots,q_n^i \in \l^3/_{\langle u\rangle},
i=1,2$ where $u \in \{x_3=0\},$ $u \neq 0.$  Suppose $S_1$ and $S_2$ are contained in $\{x_3\ge
0\}\subset\l^3/_{\langle u\rangle},$ have the same limit normal directions at the
ends and that $q_j^1=q_j^2,$ $j=1,\dots,n.$ Then $S_1 = S_2.$
\end{theorem}

\begin{proof}
The proof is based on the maximum principle. A regular maximal
surface in $\l^3$ can be represented locally as a graph $x_3 =u(x_1,x_2)$ of a smooth function $u,$ with $u_{x_1}^2 + u_{x_2}^2
< 1,$ satisfying the equation:
$$ (1-u_{x_1}^2) u_{x_2x_2} + 2u_{x_1}u_{x_2}u_{x_1x_2} + (1-u_{x_2}^2)
u_{x_1x_1}=0.$$

The maximum principle for elliptic quasilinear equations then
gives the following geometric maximum principles for maximal
surfaces:

Let $\tilde{S}_1$ and $\tilde{S}_2$ be two  maximal embedded surfaces (possibly
with boundary) in $\l^3$  which intersect tangentially at a point
$p$. Suppose that locally, around $p$, $\tilde{S}_1$ is above $\tilde{S}_2,$ that
is to say, $u_1\ge u_2$ where $u_i$ denotes the function defining
the graph $\tilde{S}_i,$ $i=1,2.$ Then $\tilde{S}_1 = \tilde{S}_2$  locally around $p$ if
one of the following hypotheses holds:
\begin{itemize}
\item $p$ is an interior point of $\tilde{S}_1$ and $\tilde{S}_2,$
\item $p$ is a boundary point of $\tilde{S}_1$ and $\tilde{S}_2$ and $\partial \tilde{S}_1$ and
$\partial \tilde{S}_2$ are tangent at $p.$
\end{itemize}

In either case, by analyticity of solutions of elliptic equations,
we also infer that the two graphs $\tilde{S}_1$ and $\tilde{S}_2$ coincide
whenever they are simultaneously defined. It is important to
emphasize that this local statement only works for maximal graphs
without singularities.\\

Consider now $S_1$ and $S_2$ as in the statement of the theorem.
Since $u$ is a horizontal vector,  vertical translations are well defined
isometries of $\l^3/_{\langle u\rangle}.$ For any $t\in \r,$ put $S_i(t)= S_i +
(0,0,t),$ $i=1,2.$

>From our assumptions and  the asymptotic behavior of the ends given in Theorem
\ref{th:representation}, we deduce that, for $t>0$ big enough, $S_1(t) > S_2,$ that
is to say $S_1(t)\cap S_2 =\emptyset$ and $S_1(t)$ is above $S_2.$ Let $t_0 = \inf
\{ t> 0\,: \, S_1(t)
> S_2\}.$ We are going to prove that $t_0 = 0.$

Suppose $t_0 >0.$
If $S_1(t_0)\cap S_2 \neq \emptyset$ then $S_1(t_0)$ and
$S_2$ have a contact point different from the singularities. But
then the interior maximum principle implies that $S_1(t_0)= S_2,$
which is absurd.

Assume now that $S_1(t_0)\cap S_2 =\emptyset$ (contact at
infinity). Because the ends of the surfaces are asymptotic to
spacelike flat half cylinders, then for $\epsilon>0$ small enough,
$S_1(t_0 -\epsilon)\cap S_2$ is a non empty compact real
1-dimensional analytic manifold containing a Jordan curve $\Gamma$
spanning two parallel annular ends without singular points $E_1
\subset S_1(t_0-\epsilon)$ and $E_2\subset S_2,$ with $E_1\cap
E_2=\Gamma.$ Let $F_1=\int_\Gamma \nu_1$ and $F_2=\int_\Gamma
\nu_2$ the flux along $\Gamma$ in $S_1(t_0-\epsilon)$ and $S_2$
resp. (see definition \ref{def:flux}). It is not hard to see that
$F_i$ is orthogonal to $u$ and to the limit normal vector at the end $E_i.$
Moreover, $\langle F_i,F_i \rangle=\langle u,u \rangle,$  and so we infer that
$F_1=F_2.$
However by the boundary maximum principle, the third coordinate of
$\nu_1$ is strictly bigger than that of $\nu_2$ along $\Gamma,$
which is a contradiction. This proves that $t_0=0,$ and reversing
the argument, that $S_1=S_2.$

\end{proof}

\begin{corollary}
The group of  ambient isometries preserving a CMF surface $S$ in
$\l^3/_{\langle u\rangle}$ coincides with:
\begin{itemize}
\item the group of orthochronous (i.e., preserving the future direction)
ambient isometries leaving invariant the set of its singularities and
preserving the set of normal directions  at the ends in case the
ends of $S$ are not parallel
\item the group of  ambient isometries leaving invariant the set of its singularities
and the limit normal vector at the ends in case $S$ has parallel
ends.
\end{itemize}
\end{corollary}

\section{The space of CMSF surfaces}
In this section we are going to study the moduli space of CMSF surfaces in $\l^3.$
By Remark \ref{re:conexion}, this space can be identified with the space
of CMF surfaces in quotients $\l^3/G,$ where $G=\langle u \rangle$ and $u$  is a spacelike vector. Moreover, as shown in Theorem \ref{th:representation},
we can restrict ourselves to the case of CMF graphs over spacelike flat cylinders in   $\l^3/_{\langle u \rangle}.$

First, we have to introduce some normalizations.

Let $\tilde{S}$ be a CMSF surface invariant by a spacelike
translation $T.$ Up to a isometry in $\l^3$ and rescaling we will
always suppose $T(p)=p+(1,0,0).$ From Theorem
\ref{th:representation} we know that $S:=\tilde{S}/_{\langle
(1,0,0) \rangle}$ is a CMF graph over the cylinder $\{x_3=0\}\subset
\l^3/_{\langle (1,0,0) \rangle}$ with flat ends. Up to a hyperbolic
rotation in $\coc$ we can suppose that one of them is asymptotic
to the half cylinder $\{x_3=0, x_2\ge 0\}.$

In the sequel we denote by $\Gg_n$ the space of CMF graphs in
$\coc$ over the cylinder $\{x_3=0\}$ with one of their ends, which will be denoted  $E_1$ in the sequel, asymptotic to
$\{x_3=0, x_2\ge 0\}$ and having $n+1$ singularities, $n\ge 1.$ We
will always suppose that all $S \in \Gg_n$ are oriented  by the
past directed normal. Note that the limit normal to the second end $E_2$
of $S$  lies in $\h^2 \cap \{x_3<0, x_1=0\}.$ The latter set is
identified, through a suitable stereographic projection, with the
real interval $]-1,+1[.$\\

Let $S\in \Gg_n$  and label $F$ as its set of singularities. By
definition, a {\em mark} in $S$ is an ordering
$\mg=(q_0,q_1,\ldots,q_n)\in \big(\coc\big)^{n+1}$ of the points
in $F,$ and we say that $(S,\mg)$ is a marked graph.
We denote by $\Mg_n$ the space of marked graphs and define the two
following maps:

$$\begin{array}{ccc}
\sg_1:\Mg_n\to\Gg_n  &\quad \mbox{and}  \quad &   \sg_2:\Mg_n\to\big(\coc\big)^{n+1}\times]-1,1[\\
 \sg_1(S,\mg)=S      & &   \sg_2(S,\mg)=(\mg,c)
\end{array}$$

where $c\in ]-1,+1[$ is the limit normal at the end $E_2$ as
explained above.\\

Label ${\mathcal P}_{n+1}$ as the symmetric group of permutations of
order $n+1$ and denote by $\mu:{\mathcal P}_{n+1} \times \Mg_n \to
\Mg_n,$  the natural action $\mu(\tau,(S,\mg)):=(S,\tau(\mg)).$
Observe that the space $\Gg_n$ can be naturally identified with
the orbit space of this action.\\


This section is devoted to prove the main result of this paper:

\begin{quote}
{\bf Main Theorem} {\em  The set $\sg_2(\Mg_n) \subset
\big(\coc\big)^{n+1} \times ]-1,+1[$ is open and
the one to one map $\sg_2:\Mg_n \to \sg_2(\Mg_n)$ provides a
global system of analytic coordinates on $\Mg_n.$

Moreover, the action $\mu$ is discontinuous and hence $\Gg_n$ has
a unique analytic structure making $\sg_1$ an analytic covering of $(n+1)!$ sheets.}
\end{quote}

This section is organized as follows: in Subsection
\ref{sub:ident} we identify $\Mg_n$ with a set
${\sb}_n\times\coc\times\{-1,1\},$ where ${\sb}_n$ is a divisor bundle associated to the Weierstrass data. The definition of
${\sb}_n$ involves some elements of classical theory of Riemann
surfaces, like the Jacobian variety and the Abel-Jacobi map, which
will be explained in this subsection.

In Subsection \ref{sub:spin}  we prove that ${\sb}_n$
has a natural structure of differentiable $(3n+1)-$manifold, and thus we use the previous identification to endow $\Mg_n$ with a
structure of differentiable manifold of dimension $3n+4.$

Finally in Subsection \ref{sub:main} we prove the Main Theorem, first
showing that $\sg_2$ is smooth when we consider the previous
differentiable structure on $\Mg_n,$ and then applying the Domain Invariance Theorem.

\subsection{Identifying $\Mg_n$}\label{sub:ident}
We split this subsection into three stages.

\subsubsection{From marked graphs in $\Mg_n$ to divisors on marked circular domains}
The following definition and notations are required.

We label ${\cal T}_n\subset \r^{3n+2}$ as the $(3n+1)$-dimensional connected analytical submanifold
consisting of points $v=(c_0,c_1,\ldots,c_n,r_0,r_1, \ldots r_n)$ in $]1,+\infty[\times\c^{n}\times(\r^+)^{n+1}$ such
that $r_0=c_0-1,$ the
discs $D_j:=\{|z-c_j| \leq r_j\},$ $j=0,1, \ldots,n,$ are pairwise disjoint and $0 \notin D_j,$ for any $j.$
We call $a_j:=\partial D_j$  and write $c_j(v):=c_j,$ $r_j(v):=r_j,$ $D_j(v):=D_j$ and $a_j(v):=a_j,$ $j=0,1,\ldots,n.$

\begin{definition}\label{def:uve}
Given $v\in {\cal T}_n,$ the domain  $\Omega(v):=\overline{\c}-\cup_{j=0}^n D_j(v)$
is said to be a marked circular domain (with $n+1$ holes).

Two marked
circular domains $\Omega(v_1)$ and $\Omega(v_2)$ are considered
equal if and only if $v_1=v_2.$
\end{definition}

For any $v\in{\cal T}_n,$ we call $\overline{\Omega(v)}^*$ and $N(v)$ the
mirror of $\overline{\Omega(v)}:=\Omega(v)\cup \big(\cup _{j=0}^n
\partial a_j(v)\big)$ and  the double surface of $\overline{\Omega(v)},$ respectively. Recall that
$$N(v)=\overline{\Omega(v)} \cup \overline{\Omega(v)^*}$$
with the identification
$\partial(\overline{\Omega(v)})\equiv\partial(\overline{\Omega(v)^*}).$
We know that $N(v)$ is a closed Riemann surface, and
$\overline{\Omega(v)^*} \cap \overline{\Omega(v)}$
consists of the $n+1$ analytic circles $a_j(v),$ $j=0,\ldots,n.$
Moreover, we denote by $J_v:N(v) \to N(v)$ the antiholomorphic
involution applying any point to its mirror image. Note that the
fixed point set of $J_v$ coincides with $\cup_{j=0}^n a_j(v).$

\begin{figure}[hbtp]
\begin{center}
\includegraphics[width=14cm]{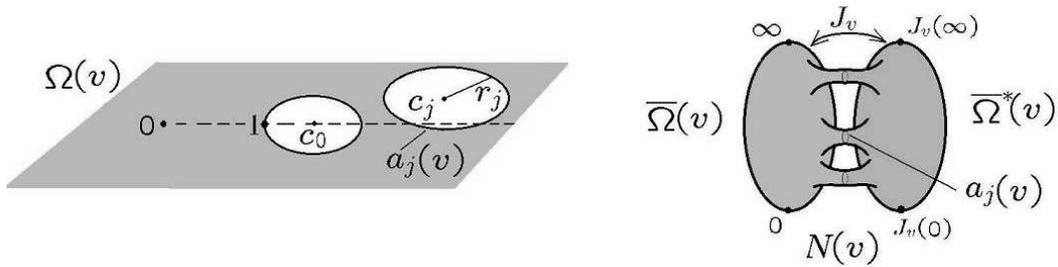}
\caption{$\Omega(v),$ $N(v)$ and $J_v.$ \label{fig:omega}}
\end{center}
\end{figure}

\begin{remark} \label{re:mirror}
A conformal model for $\overline{\Omega(v)}^*,$ $v \in {{\cal
T}_n},$ consists of the planar domain
$\overline{\Omega(v)}^*:=\{J_{v}(z) \;:\; z \in
\overline{\Omega(v)}\},$ where
$$J_{v}(z):=c_0(v) +\frac{r_0(v)^2}{\bar{z}-c_0(v)}$$
is the Schwarz reflection about $a_0(v)=\{|z-c_0(v)|=r_0\}.$
Moreover, $N(v)$ can be identified to the quotient of
$\overline{\Omega(v)} \cup \overline{\Omega(v)}^*$ under the
identification $z \equiv J_{v}(z),$ $z \in \partial \Omega(v).$
\end{remark}

Let $Y=(S,\mg) \in \Mg_n.$ By Theorem \ref{th:representation} we
know that the conformal support of $S$ with the prescribed
orientation is biholomorphic to a twice punctured circular domain,
where the two punctures $\{0,\infty\}$ correspond to the ends and
the boundary circles to the singularities (without loss of generality, $E_1$ corresponds to $z=0$).

We can therefore
associate to $(Y,\mg)$ a unique element $\Omega(v) \in {\cal T}_n$ and conformal
immersion $X: \overline{\Omega(v)}-\{0,\infty\} \to \coc,$ as stated in the following lemma:

\begin{lemma}\label{lem:ident}
Given $Y=(S,\mg) \in \Mg_n,$ where $\mg=(q_0\ldots,q_n),$  there
are {\em unique} $v\in {\cal T}_n$ and conformal maximal immersion
$X:\overline{\Omega(v)}-\{0,\infty\} \to \coc$ such that:
\begin{enumerate}[(i)]
\item
$S - F$ is biholomorphic to $\Omega(v)-\{0,\infty\}$ (in the
sequel, they will be identified),
\item
$S=X(\overline{\Omega(v)}-\{0,\infty\}),$
\item
$z=0$ correspond to the end $E_1,$
\item
$q_j=X(a_j(v)),$  $j=0,\ldots,n.$
\end{enumerate}
\end{lemma}

Lemma \ref{lem:sing} and Theorem \ref{th:representation}
also give that the Weierstrass data of $X,$ $(g,\phi_3),$ satisfy the
symmetries $g\circ J_v=1/\overline{g}$ and
$J_v^*(\phi_3)=-\overline{\phi}_3,$ and that $g$ has exactly $n+1$
zeros $0,w_1,\ldots,w_n \in \Omega(v)$ counted with multiplicity.

Therefore, putting $D=w_1\cdot\ldots\cdot w_n\in Div_n(\Omega(v)),$ 
it is easy to see that divisors for the Weierstrass data must be:
\begin{equation}\label{eq:divisors}
(g)=\frac{D\cdot 0}{J_v(D\cdot 0)}\quad\mbox{and}\quad
(\phi_3)=\frac{D\cdot J_v(D)}{\infty\cdot  J_v(\infty)}
\end{equation}

Since the divisor $D$  determines uniquely the data
$(g,\phi_3)$ up to multiplicative constants, and these data control the immersion $X,$
we infer that the couple $(v,D)$ encloses all the information about the surface.

\subsubsection{The Abel-Jacobi map on the  bundle of divisors}

In order to understand the moduli space $\Mg_n,$ it is crucial to control the structure
of the family of couples $(v,D)$ for
which there exist Weierstrass data $(g,\phi_3)$ satisfying Equation (\ref{eq:divisors}).
The Abel-Jacobi map (defined below) will play here  a fundamental role.

We need some extra
notation.


Given a Riemann surface $\rb,$ we denote by
$$Div_k(\rb)=\{D\;:\; D\, \mbox{ is an integral multiplicative divisor on}\, \rb\, \mbox{of degree}\, k \}$$
Recall that $Div_k(\rb)$ is
the quotient of $\rb^k$ under the action of the group of
permutations of order $k,$ and we denote by $p_k:\rb^k \to Div_k(\rb)$
the canonical projection. We endow $Div_k(\rb)$ with the natural
analytic structure induced by $p_k.$

In what follows, for any $k\in\n,$ we denote by
$$\div_k=\bigcup_{v\in{\cal T}_n} Div_k(\Omega(v))=\{(v,D)\;:\;v\in{\cal T}_n,\,D\in Div_k(\Omega(v))\}$$
and we refer to it as the {\em bundle of k-divisors}.

Obviously, $\div_k$ is a real analytical manifold
(see \cite{fer-lop-sou} for more details).\\

\begin{figure}[hbtp]
\begin{center}
\includegraphics[width =7.5cm]{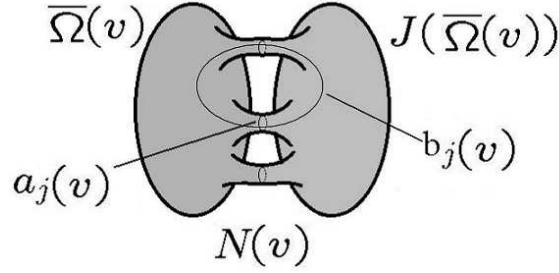}
\caption{The canonical homology basis $B(v)$
\label{fig:homologia}}
\end{center}
\end{figure}

Let ${\J}(v)$  be the Jacobian variety of the compact Riemann
surface $N(v)$ associated to the following canonical homology
basis:

We identify the homology classes of the boundary circles $a_j(v)$
in $\Omega(v)$ with their representing curves $j=0,1,\ldots,n.$
Note first that $J_v$ fixes $a_j(v)$ pointwise, and so,
$J_v(a_j(v))=a_j(v).$ Take a curve $\gamma_j \subset
\overline{\Omega(v)}$ joining $a_0(v)$ to $a_j(v),$ in such a way
that the curve $b_j(v)$ obtained by joining $\gamma_j$ and
$J_v(\gamma_j)$  satisfies that the intersection numbers
$(b_j(v),b_h(v))$ vanish, and $(a_j(v),b_h(v))=\delta_{jh},$ where
$\delta_{jh}$ refers to the Kronecker symbol. Observe that
$J_v(b_j(v))=-b_j(v)$ in the homological sense, and its homology
class does not depend on the choice of $\gamma_j.$ In other words,
the identity $J_v(b_j(v))=-b_j(v)$ characterizes
$B(v)=\{a_1(v),\ldots, a_n(v),b_1(v),\ldots,b_n(v)\}$ as canonical
homology basis of $N(v)$ (see Figure \ref{fig:homologia}).

Call $\{\eta_1(v),\ldots,\eta_n(v)\}$ the dual basis of $B(v)$ for
the space of holomorphic 1-forms on $N(v),$ that is to say, the
unique basis satisfying $\int_{a_k(v)}\eta_j(v)=\delta_{jk},$
$j,k=1,\ldots,n,$ and put $\Pi(v)=(\pi_{j,k}(v))_{j,k=1,\ldots,n}$
for the associated matrix of periods, $\pi_{j,k}(v)=\int_{b_j(v)}
\eta_k(v).$

Then the Jacobian variety of $N(v)$ is ${\J}(v)=\c^n/L(v),$ where
 $L(v)$ is the lattice over $\z$
generated by $\{e^1,\ldots,e^n,\pi^1(v),\ldots,\pi^n(v)\},$ where
$$e^j=\,^T(0,\ldots, \stackrel{j}{1},\ldots,0)\quad
\mbox{and}\quad
\pi^j(v)=\,^T(\pi_{1,j}(v), \ldots,\pi_{n,j}(v))$$
The {\em
Jacobian bundle} is defined as
$${\J}_n=\bigcup_{v\in{\cal T}_n}{\J}(v)$$
${\J}_n$ has a natural structure of analytic manifold (see \cite{fer-lop-sou} for more details).\\

For any $v \in {\cal T}_n,$ we call $\varphi_v:N(v) \to {\J}(v)$
the Abel-Jacobi embedding defined by
$$\varphi_v(z)=\pg_v\left( \int_{1}^{z} {^T(\eta_1(v),\ldots,\eta_n(v))}\right)$$
where $\pg_v:\c^n\to{\J}(v)$ is the canonical projection (recall
that $1\in \overline{\Omega(v)} \subset N(v)$ uniformly on $v$).
We extend $\varphi_v$ with the same name to the Abel-Jacobi map
$\varphi_v:Div_k(N(v)) \to {\J}(v)$ given by
$$\varphi_v(P_1 \cdot
\ldots \cdot P_k)= \sum_{j=1}^k \varphi_v(P_j),\qquad k \geq 1.$$

We also define $\varphi:\div_k\to{\J}_n$ by $\varphi(v,D)=(v,\varphi_v(D))$

\subsubsection{The identification $\Mg_n \equiv {\cal
S}_n\times\coc\times\{-1,1\}$}

Summarizing, we know that given $Y=(S,\mg)\in\Mg_n$ its associated
Weierstrass data (defined on $N(v)$ for the unique $v\in{\cal
T}_n$ given in Lemma \ref{lem:ident} ) satisfy Equation
(\ref{eq:divisors}). Abel's Theorem gives
$$\varphi_v(D \cdot 0)-\varphi_v(J_v(D \cdot \infty))=0\quad\mbox{and}\quad
\varphi_v(D\cdot J_v(D))-\varphi_v(\infty \cdot
J_v(\infty))=T(v)$$
where $T(v)\in{\J}(v)$ is the image by
$\varphi_v$ of the divisor associated to a meromorphic 1-form on
$N(v).$ By Abel's theorem, $T(v)$ is independent of the choice of
the meromorphic 1-form (see \cite{farkas}).

These two equations lead to $2 \varphi_v(D\cdot
0)-\varphi_v(0\cdot\infty\cdot J_v(0)\cdot J_v(\infty))=T(v) $

Therefore, it is natural to define, for any $v\in{\cal T}_n,$
$${\cal S}_n(v):=\{D\in Div_n(\Omega(v))\;:\;
2 \varphi_v(D\cdot 0)=T(v) + \varphi_v(0\cdot\infty\cdot
J_v(0)\cdot J_v(\infty))\},$$ and also
$${\cal S}_n:=\{(v,D)\;:\;v\in{\cal
T}_n,\; D\in{\cal S}_n(v)\}$$
We refer to ${\cal S}_n$ as the {\em spinorial bundle.}\\

\begin{definition} \label{def:cale}
With the previous notation, we call $\E$ the map given by
$${\E}: \Mg_n \to {\cal S}_n\times \coc\times\{-1,1\} $$
$${\E}(S,\mg)= \big( (v,D),q_0,\varepsilon_0 \big)$$
where $\mg=(q_0,\ldots,q_n)$ and
$\varepsilon_0\in\{-1,1\}$ is the sign of the third coordinate of the flux at $q_0$
(see Definition \ref{def:flux})
\end{definition}

The main goal of this paragraph is to show that ${\cal E}$ is bijective. Note that the first coordinate of  ${\cal E}$
encloses the information about the conformal structure and
Weierstrass data of the marked graph, while the second one is
simply translational. The third coordinate has been introduced just for distinguishing between
two graphs having $q_0=0$ and being symmetric  with respect to the plane $\{x_3=0\}.$

The following notation and lemmae are required. \\

Consider the holomorphic 1-form $\overline{J_v^*(\eta_j(v))}.$
Taking into account that $J_v$ fixes $a_j(v)$ pointwise, we infer
that $\int_{a_k(v)} \overline{J_v^*(\eta_j(v))} =\delta_{jk},$ and
so, $J_v^*(\eta_j(v))=\overline{\eta_j(v)}.$  Moreover, since
$J_v(b_j(v))=-b_j(v),$ then $\pi_{j,k}(v)=\int_{b_k} \eta_j(v)$ is
an imaginary number, for any $j$ and $k.$

It follows that there exists a unique analytic mirror involution
$I_v:{\J}(v) \to {\J}(v)$ satisfying $I_v (\pg_v(w))= \pg_v
(\overline{w}),$ for any $w \in \c^n.$ Moreover, as $J_v(1)=1$  then
$\varphi_v\circ J_v= I_v\circ\varphi_v.$

We call ${\cal I}:{\cal J}_n\to{\cal J}_n,$ the map given by
${\cal I}(v,\pg_v(w))=(v,I_v(\pg_v(w))).$

\begin{lemma}[\cite{fer-lop-sou}]\label{lem:todas}
The  maps $\varphi:\div_k\to{\cal J}_n,$
$\varphi(v,D)=(v,\varphi_v(D)),$ and $\hat{T}:{{\cal T}}_n\to{\cal
J}_n,$ $\hat{T}(v)=(v,T(v))$ are smooth.
\end{lemma}

As a consequence of the smoothness  of $\hat{T}$ and $\varphi$ it
follows that there are exactly $2^{2n}$ differentiable maps
$\hat{E}_1,\ldots,\hat{E}_{2^{2n}}:{\cal T}_n\to{\cal J}_n,$
$\hat{E}_j(v)=(v,E_j(v)),$ satisfying $2E_j(v)=T(v)+\varphi(v,0\cdot\infty\cdot J_v(0)\cdot J_v(\infty))$ for any
$j.$

The next result shows that these {\em spinor sections} are
invariant under the mirror involution. This fact will be crucial for recovering
the Weierstrass data from an element in the spinorial bundle.

\begin{lemma} \label{lem:spinor}
${\cal I} \circ \hat{E}_j=\hat{E}_j,$ for any $j=1,\ldots,2^{2
n}.$

\end{lemma}
\begin{proof}
Indeed, note that $I_v(E_j(v))=E_j(v)+\pg_v\big(\frac{1}{2}
\sum_{h=1}^n (m_h(v) e^h+n_h(v) \pi^h(v))\big),$ where $m_h(v),$
$n_h(v) \in \z$ are continuous functions of $v.$ Using that
${{\cal T}_n}$ is connected we get that $m_h(v),\;n_h(v)$ are
constant. Hence, the set ${\cal A}_j:= \{v \in {{\cal T}_n} \;:\;
I_v(E_j(v))=E_j(v)\}$ is either empty or the whole of ${{\cal
T}_n}.$ On the other hand, $E_j(v)=E_1(v)+q_j(v),$ where $2
q_j(v)=0,$ and so, $I_v(q_j(v))=q_j(v).$  Therefore ${\cal
A}_1={{\cal T}_n}$ if and only if ${\cal A}_j={{\cal T}_n}$ for
any $j.$

Consider the compact genus $n$ Riemann surface $N=\{(z,w) \in
\overline{\c}\;:\; w^2=\prod_{i=1}^{2n+2} (z-c_i)\},$ where $c_i
\in \r$ and $c_1<c_2< \ldots <c_{2n}<0<c_{2n+1}<c_{2n+2}.$ 
The function $w$ has a well defined branch $w_+$ on the planar
domain $\Sigma =\overline{\c}- \cup_{i=0}^{n} [c_{2
i+1},c_{2i+2}].$ Moreover there exists a biholomorphism from the
domain $\{(z,w_+(z)) \;:\; z \in \Sigma\} \subset N$  to a
circular domain $\Omega(v_0),$ $v_0 \in {{\cal T}_n}$ taking
$0_+:=(0,w_+(0))$ to $0$ and $\infty_+:=(\infty,w_+(\infty))$ to
$\infty.$ Furthermore, up to this biholomorphism, $N=N(v_0)$ and
$J=J_{v_0}$ is given by $J(z,w)=(\overline{z},-\overline{w}).$

Define the meromorphic 1-form $\nu=\prod_{i=1}^{n+1} (z-c_i)
\frac{dz}{zw}$ on $N(v_0)$ and observe that its canonical divisor
is given by
$$(\nu)=\frac{c_1^2 \cdot \ldots \cdot
c_{n+1}^2}{0_+\cdot J_{v_0}(0_+)\cdot \infty_+\cdot
J_{v_0}(\infty_+)}$$
where we are identifying $c_i \equiv (c_i,0)
\in N(v_0).$  Since $J_{v_0}(c_i)=c_i,$ then
$l_0:=\sum_{i=1}^{n+1} \varphi_{v_0} (c_i) \in {\cal J}({v_0})$ is
invariant under $I_{v_0}$ and  $2
l_0=T(v_0)+\varphi_{v_0}(0_+\cdot J{v_0}(0_+)\cdot \infty_+\cdot
J_{v_0}(\infty_+)).$ Up to relabelling, we can suppose that
$l_0=E_1(v_0)$ and hence ${\cal A}_1={{\cal T}_n}.$ This completes
the proof.
\end{proof}

\begin{proposition}\label{pro:cale}
The map ${\cal E}: \Mg_n \to {\cal S}_n\times\coc\times\{-1,1\}$ is bijective.

\end{proposition}

\begin{proof}

If $x \in {\cal S}_n,$ $x=(v,D),$ then $\varphi_v(D\cdot
0)=E_i(v)$ for some $i\in\{1,\ldots,2^{2n}\}.$

Since $I_v(E_i(v))=E_i(v),$ (Lemma \ref{lem:spinor}) we have
$\varphi_v(D\cdot 0)-\varphi_v(J_v(D\cdot 0))=0.$ By Abel's
theorem, there exists a unique meromorphic function $g_x^0$ of
degree $n+1$ on $N(v)$ satisfying
\begin{equation}\label{eq:g}
(g_x^0)=\frac{D \cdot 0}{J_v(D \cdot 0)}\quad \mbox{and}\quad
g_x^0(1)=1
\end{equation}
Observe that since $J_v(1)=1$ we have $g_x^0\circ
J_v=1/\overline{g_x^0}.$

On the other hand, as $\varphi_v(D \cdot
J_v(D))-\varphi_v(J_v(\infty) \cdot \infty)=T(v),$ then there
exists a meromorphic 1-form $\phi$ on $N(v)$ with canonical
divisor $\frac{D \cdot J_v(D)}{\infty \cdot J_v(\infty)}.$ Up to a
multiplicative constant we can suppose that $\phi$ satisfies
$J_v^*(\phi)=-\overline{\phi}.$ If we write $\phi(z)=h(z)\,
\frac{dz}{z-c_0(v)},$ $z\in
U(v)=\big(\Omega(v)-\{0,\infty\}\big)\cup\big(\Omega(v)^*-\{J_v(0),J_v(\infty)\}\big)\cup
a_0(v),$ we infer that $h(z) \in \r^*,$ for any $z$ satisfying
$|z-c_0(v)|=r_0(v)=c_0(v)-1.$ Then, define
$$\phi_3^0(x):=\frac{1}{h(1)}\, \phi,$$ and observe that the
equations
\begin{equation}\label{eq:phi}
(\phi_3^0(x))=\frac{D \cdot J_v(D)}{\infty \cdot
J_v(\infty)}\quad\mbox{and}\quad h_3^0(1)=1
\end{equation}
characterize $\phi_3^0(x)$ as meromorphic 1-form on $N(v).$

Given $\chi=\big( x,q_0,\varepsilon_0 \big) \in {\cal S}_n\times
\coc\times\{-1,1\}$ it is easy to check that there exist unique
$\theta_\chi\in \{|z|=1\}$ and $r_\chi\in \r$ such that, for $\phi_3(\chi)=r_\chi\phi_3^0(x)$ and $
g_\chi=\theta_\chi g_x^0,$ the map
$X_\chi:\overline{\Omega(v)} \to \coc,$
$$X_\chi(z):= q_0+\mbox{Real} \int_1^z
\big(\frac{i}{2}(\frac{1}{g_\chi}- g_\chi),-\frac{1}{2}(\frac{1}{g_\chi}+
g_\chi), 1\big)\phi_3(\chi) $$
is well defined and determines a CMF graph
$S_\chi:=X_\chi\Big(\overline{\Omega(v)}-\{0,\infty\}\Big)
\in\Gg_n$ satisfying $X(a_0(v))=q_0$ and $\epsilon_0$ is the sign of the third coordinate
of the flux along $a_0(v).$

Defining the mark $\mg_\chi$ by $q_{j}=X_\chi(a_j(v)),$
$i=0,\ldots,n,$ it is now clear that ${\cal
E}^{-1}(\chi)=\{(S_\chi,\mg_\chi)\},$ and so, ${\cal E}$ is
bijective.
\end{proof}

\subsection{Structure of the spinorial bundle ${\sb}_n$}\label{sub:spin}

In the previous subsection we have identified $\Mg_n$ with the
space ${\cal S}_n\times\coc\times\{-1,1\}.$ Our aim now is to show
that ${\cal S}_n,$ and so $\Mg_n,$ has a natural structure of differentiable
manifold.

\begin{theorem}[Structure of the spinorial bundle] \label{th:submanifold}
The space ${\cal S}_n$ is a smooth real $(3n+1)$-dimensional
submanifold of $\div_{n}$ and the map $\vg:{\cal S}_n \to {\cal
T}_n,$ $\vg(v,D) = v$ is a finite covering.
\end{theorem}

\begin{proof}
The fact ${\cal S}_n \neq \emptyset$ follows from the existence of
CMSF surfaces with an arbitrary number of singularities  in the
quotient (see \cite{fer-lop}) and Proposition \ref{pro:cale}. The
key step of this proof is that ${\cal S}_n$ does not contain any
{\em special} divisor (see \cite{farkas}).

Consider the differentiable map $H: \div_{n} \to {\cal J}_n$ given
by
$$H(v,D)=\Big(v,2\varphi_v(D)-\varphi_v\big(0 \cdot J_v(0)\cdot \infty
\cdot J_v(\infty)\big)-T(v)\Big),$$
and note that ${\cal
S}_n=\{(v,D) \in \div_{n} \;:\;H(v,D)=(v,0)\}.$ In order to prove
that ${\cal S}_n$ is a differentiable submanifold of $\div_{n},$
it suffices to check that $dH_q$ is bijective at any point $q$ of
${\cal S}_n.$

Let $q_0:=(v_0,D_0)$ be an arbitrary point of ${\cal S}_n.$
Observe that $dH_{q_0}$ is bijective if and only if the map
 $$H_0:Div_n(\Omega(v_0))\to {\cal J}(v_0)$$
 $$D\mapsto \varphi_{v_0}(D)-\varphi_{v_0}(D_0) $$
is a local diffeomorphism at $D_0.$

We are going to write the expression of $H_0$ in local coordinates
around $D_0 \in Div_n(\Omega(v_0))$ and $H_0(D_0)=0 \in {\cal
J}(v_0).$ To do this, write $D_0=z_1^{m_1} \cdot\ldots \cdot
z_s^{m_s} \in Div_n(\Omega(v_0)),$ and denote by
$(U_j,w_j:=z-z_j)$ the conformal parameter in $\Omega(v_0),$ where
$U_j$ is the open disc of radius $\epsilon>0$ centered at $z_j,$
$j=1,\ldots,s.$ Put $U=\prod _{j=1}^s
U_j^{m_j}\subset\Omega(v_0)^n=\Omega(v_0)\times\stackrel{n}{\ldots}\times\Omega(v_0).$
Then $p_n(U)$ is a neighborhood of $D_0$ in $Div_n(\Omega(v_0)),$
where $p_n:\Omega(v_0)^n\to Div_n(\Omega(v_0))$ is the projection
associated to the action of the group of permutations of order $n$
on $\Omega(v_0)^n.$ A coordinate chart for $Div_n(\Omega(v_0))$
around $D_0$ is given by
$$ \xi:p_n(U)\to\c^n, \qquad  \xi (\prod_{j=1}^s Q_{1,m_j} \cdot \ldots \cdot
Q_{m_j,m_j})=((t_{1,m_j},\ldots,t_{m_j,m_j})_{j=1, \ldots,s}),$$
where $t_{h,m_j}=\sum_{l=1}^{m_j} (z_j(Q_{l,m_j}))^h,$
$h=1,\ldots,m_j,$ $j=1,\ldots,s.$ For more details see
\cite{farkas}.

Label $\pg:\c^n\to{\cal J}(v_0)=\c^n/L(v_0)$ as the natural
projection and consider a neighborhood $W'$ of $H(D_0)=0$ such
that $\pg\, :\, W:=\pg^{-1}(W)\to W'$ is a diffeomorphism and
$H_0(U)\subset W'.$

Write $\eta_i(v_0)(w_j)=f_{i,j}(w_j)dw_j$ on $W_j:=w_j(U_j)$ for
$i=1,\ldots,n,$ $j=1,\ldots s.$ The local expression $\hat{H}_0$
of $H_0$ around $D_0,$ $\hat{H}_0=\pg^{-1}\circ H_0\circ
\xi^{-1},$ is given by
$$\hat{H}_0:\xi(p_n(U))  \to W$$
$$\hat{H}_0(t)= \sum_{j=1}^s \sum_{h=1}^{m_j} \int_{0}^{w_{h,m_j}}
f_{j}(w_j)dw_j $$
where $f_{j}=\,^T(f_{1,j},\ldots,f_{n,j}),$
$w_{h,m_j}\equiv w_j,$ $t_{l,m_j}=\sum_{h=1}^{m_j}w_{h,m_j}^l,$
$h=1,\ldots,m_j,$
$t=(t_{1,m_j},\ldots,t_{m_j,m_j})_{j=1,\ldots,s}.$

Put $f_j(w_j)=\sum_{l=0}^{\infty} b_{j,l} w_j^l,$
$b_{j,l}\in\c^n,$ $j=1\ldots s.$ Then the Taylor series for the
holomorphic map $w_{h,m_j}\mapsto
\int_{0}^{w_{h,m_j}}f_{j}(w_j)dw_j$ is
$\int_{0}^{w_{h,m_j}}f_{j}(w_j)dw_j=\sum_{l=1}^\infty a_{j,l}
w_{h,m_j}^l,$ where $a_{j,l}=\frac{1}{l} b_{j,l-1},$ $l\geq 1,$
$j=1,\ldots,s.$ It is not hard to check that $\hat{H}_0(t)\sum_{j=1}^s \sum_{l=1}^{m_j} a_{j,l} t_{l,m_j}+R(t),$ where the
first derivatives of $R$ with respect to $t_{l,m_j}$ vanish at
$t=0,$ and so the column vectors of the Jacobian matrix of
$\hat{H}_{0}$ are $\{a_{l,j},l=1,\ldots,m_j,j=1\ldots,s\}.$

Reasoning by contradiction, suppose that the rows of that matrix
are linearly dependent, which is equivalent to saying that there
exists a holomorphic 1-form $\omega_0$ in $N(v_0)$ having a zero
at $z_j \in \Omega(v_0) \subset N(v_0)$ of order at least $m_j,$
$j=1,\ldots,s.$ A direct application of Riemann-Roch theorem gives
the existence of a non-constant meromorphic function $f$ on
$N(v_0)$ having poles at $z_1,\ldots,z_s$ with order at most
$m_1,\ldots,m_s,$ respectively. In particular, $f$ has degree less
than or equal to $n.$ As $J_{v_0}(0)$ is not a pole of $f,$ up to
adding a constant we can suppose that
  $f(J_{v_0}(0))=0.$

On the other hand, since
$\varphi_{v_0}(D_0)+\varphi_{v_0}(0)=E_i(v_0)$ and
$I_{v_0}(E_i(v_0))=E_i(v_0),$ $i \in \{1,\ldots,2^{2n}\},$ then we
get $\varphi_{v_0}(D_0 \cdot 0)-\varphi_{v_0}(J_{v_0}(D_0 \cdot
0))=0.$ Therefore, a direct application of Abel's theorem  gives
the existence of a meromorphic function $g$ of degree $n+1$ on
$N(v_0)$ whose principal divisor coincides with $\frac{D_0 \cdot
0}{J_{v_0}(D_0 \cdot 0)}.$ As $J_{v_0}$ is an antiholomorphic
involution with fixed points, it is not hard to check that $g
\circ J_{v_0}=r/\overline{g},$ $r>0.$ Hence, up to multiplying $g$
by the factor $r^{-1/2},$ we can suppose that $g \circ
J_{v_0}=1/\overline{g}.$

Note that $g \in {\cal F}_{v_0},$ where for any
 $v \in {\cal T}_n,$ ${\cal F}_v$ denotes the family
of meromorphic functions $h$ of degree $n+1$ in $N(v)$ with zeroes
in  $\Omega (v)$ and satisfying $h\circ J_{v}=1/\overline{h}.$

\vspace{0.3cm}
\begin{quote}{{\bf Claim:}  \em Let  $f_\lambda $ denote the meromorphic
function $\frac{1+\lambda f}{1+\overline{\lambda (f \circ
J_{v_0})}},$ $\lambda \in \c.$ Then, $f_\lambda$ is not constant,
for any $\lambda \in \c^*.$ Moreover, $g_\lambda:=g f_\lambda \in
{\cal F}_{v_0}$ for any $\lambda \in \c.$}
\end{quote}

Assume $f_\lambda=c,$ where  $c,$ $\lambda \in \c^*.$ Then, we
infer that $1+\lambda f=c(1+\overline{\lambda (f\circ J_{v_0})})$
and so the polar divisor of $f,$ which is contained in $D_0,$ is
invariant under $J_{v_0}.$ This is absurd because $D_0 \in
Div_{n}(\Omega(v_0))$ and $\Omega(v_0) \cap
J_{v_0}(\Omega(v_0))=\emptyset.$

For the second part of the claim, first note that the principal
divisor of $g_{\lambda}$ is
 $(g_{\lambda})= \frac{D_{\lambda}\cdot\infty}
{J_{v_0}(D_{\lambda})\cdot J_{v_0}(\infty)},$ where $D_{\lambda}$
is
 an integral divisor of degree $\leq n$ and so
 the degree of $g_\lambda$
is $\leq n+1,$ $\lambda \in \c.$  Moreover,  $g_{\lambda}$ is not
constant for any $\lambda$ (otherwise, $J_{v_0}(0)$ would be
 a zero of $1+\lambda f,$ contradicting  $f(J_{v_0}(0))=0$).

Let $A$ be the set $\{\lambda \in \c \;:\; g_\lambda \in {\cal
F}_{v_0}\},$ and observe that $0\in A.$ It suffices to see that
$A$ is open and closed.

The openness of $A$ is  an elementary consequence of Hurwitz
theorem (we are using the fact  that the degree of $g_{\lambda}$
is at most $n+1$).

Finally, let us prove  that $A$ is closed. Let $\lambda_0 \in
\overline{A},$ and take $\{\lambda_n\}_{n \in \n} \to \lambda_0,$
where $\{\lambda_n \;:\; n \in \n\} \subset A.$ The sequence
$\{g_n:=g_{\lambda_n} \}_{n \in \n}$ converges to
$g_0:=g_{\lambda_0}$ uniformly on $N(v_0).$ We know  that  $g_n
\circ J_{v_0}=1/\overline{g_n}$ and so the zeros of $g_n$
 lie in
$\Omega(v_0),$ therefore, $g_n$ is holomorphic on $\Omega(v_0),$
$n\in\n$ and so the same holds for $g_0.$ Moreover, since
$|g_0|=1$ on $\partial\Omega(v_0)$ and it is non constant, the
maximum principle implies that  $|g_0|<1$  on $\Omega(v_0)$ and we
infer that $g_0$ has no critical points on $\partial
\Omega_{v_0}.$  As $\partial \Omega_{v_0}$ consists of $n+1$
disjoint circles, this means that $g_0$ takes on any complex
number $\theta \in \s^1$ at least $n+1$ times. Hence the degree of
$g_0$ must be $n+1$ and $g_0 \in {\cal F}_{v_0}.$ This concludes
the proof of the claim.

\vspace{0.3cm}

To get the desired contradiction  take $P \in \partial
\Omega(v_0)$ such that $f(P) \neq 0,\;\infty,$  and choose
$\lambda'= \frac{-1}{f(P)}.$ Since $J_{v_0}(P)=P,$ the meromorphic
function $g_{\lambda'}$ has degree less than $n+1,$ and so,
$\lambda' \notin A=\c,$ which is absurd.

Summarizing, we have proved that ${H}|_{{\cal S}_n}: {\cal S}_n
\to {\bf 0},$ ${H}(v,D)=(v,0),$ is a local diffeomorphism, where
${\bf 0}=\{(v,0)\;:\;v\in {\cal T}_n \}\subset{\cal J}_n$ is the null
section in the Jacobian bundle.  Consequently, the projection
$\vg:{\cal S}_n \to {\cal T}_n,$ $\vg(v,D)=v,$ is a local
diffeomorphism too. To finish, it suffices to check that  $\vg$ is
also proper. Indeed, take a sequence
$\{(v_k,D_k)\}_{k\in\n}\subset {\cal S}_n$ such that
$\{v_k\}_{k\in\n}$ converges to a point $v_{\infty}\in {\cal
T}_n.$ We can assume that $\varphi_{v_k}(D_k\cdot 0)=E_i(v_k)$ for
any $k\in\n.$  Since $I_{v_k}(E_i(v_k))=E_i(v_k),$ we get
$\varphi_{v_k}(D_k \cdot 0)-\varphi_{v_k}(J_{v_k}(D_k\cdot 0))
=0.$ By Abel's theorem  there is a meromorphic function $g_{k} \in
{\cal F}_{v_k}$ with canonical divisor $\frac{D_k \cdot
0}{J_{v_k}(D_k\cdot 0)}.$

Let us see that, up to taking a subsequence, $\{g_k\}_{k \in \n}
\to g_\infty \in {\cal F}_{v_\infty}.$ Reflecting about all the
components of $\partial \Omega(v_k),$  we can meromorphically
extend $g_k$  to a {\em planar} open neighborhood $W_k$ of
$\overline{\Omega(v_k))},$ $k \in \n.$ By continuity and for $k_0$
large enough, the set $W=\cap_{k \geq k_0} W_k$ is a planar
neighborhood of $\overline{\Omega(v_\infty)}.$ Classical normality
criteria show that, up to taking a subsequence, $\{g_k\}_{k\in\n}$
converges uniformly on $\overline{\Omega(v_\infty)}$ to a function
$g_{\infty}$ which is meromorphic beyond
$\overline{\Omega(v_\infty)}.$  It is clear that $|g_{\infty}|= 1$
on $\partial\Omega(v_{\infty}),$ $|g_{\infty}|< 1$ on
$\Omega(v_{\infty})$ and $g_{\infty}(0)=0.$ This proves that
$g_{\infty}$ is non constant and can be extended to $N(v_\infty)$
by the Schwarz reflection $g_{\infty}\circ
J_{\infty}=1/\overline{g_{\infty}}.$ Since $\deg(g_k)=n+1,$ then
Hurwitz theorem implies that $\deg(g_{\infty})\leq n+1.$ On the
other hand, $|g_{\infty}|=1$ only on $\partial
\Omega(v_{\infty}),$ and so $g_\infty$ is injective on every
boundary component of $\Omega(v_\infty).$ Therefore, the degree of
$g_\infty$ must be exactly $n+1$ and $g_\infty \in {\cal
F}_{v_\infty}.$

Finally, note that $\frac{D_{\infty} \cdot 0}{
J_{v_{\infty}}(D_{\infty}\cdot 0)}$ where $D_\infty \in \div_n$
and use Hurwitz theorem to infer that  $\{D_k\}_{k\in\n} \to
D_{\infty}\in\div_n.$ Since ${\cal S}_n$ is a closed subset of
$\div_n,$ we get $D_{\infty}\in{\cal S}_n(i),$ which proves the
properness of $\vg:{\cal S}_n \to {\cal T}_n$ and so the theorem.

\end{proof}

\subsection{Proof of the Main Theorem}\label{sub:main}

In the preceding section, we have endowed $\Mg_n$ of a differentiable structure.
It is natural to ask whether $\sg_2:\Mg_n\to\r^{3n+4}$ is a smooth map or not.

In order to do this, we have to show that the Weierstrass data of  elements in
$\Mg_n$ depends smoothly on its associated divisor in ${\cal
S}_n.$ This requires a convenient concept of
differentiability for maps from the bundle of divisors to the space
of meromorphic functions or 1-forms. The first part of this
subsection is devoted to present these concepts.\\

For any $v \in {\cal T}_n,$ call ${\cal M}(v)$ the family of
meromorphic functions on $N(v).$ The corresponding  bundle over
${\cal T}_n$ is denoted by ${\cal M}_n=\cup_{v \in {\cal T}_n}
{\cal M}(v).$

Likewise, we call  ${\cal H}(v)$  the space of  meromorphic
$1$-forms on $N(v)$ and denote  by ${\cal H}_n=\cup_{v \in {{\cal
T}_n}} {\cal H}(v)$ the associated bundle over ${\cal T}_n.$
\\

Given $v\in{\cal T}_n$ and $k_1,$ $k_2\in\n$ we denote by
$Div_{k_1,k_2}(v)$ the product manifold,
$$Div_{k_1,k_2}(v)=Div_{k_1}(v)\times Div_{k_2}(v),$$
and by $\div_{k_1,k_2}$ its associated bundle over ${\cal T}_n,$
$\div_{k_1,k_2}=\cup_{v\in{\cal T}_n} Div_{k_1,k_2}(v).$ Like in
the case of $\div_k,$ $\div_{k_1,k_2}$ has a natural structure of
analytical manifold. We use the convention $\div_{k,0}=\div_k$ and
$\div_{0,0}={\cal T}_n.$

\begin{definition}[smoothness with $k$-regularity]\label{def:kreg}
Let $M_j$ be a real manifold of dimension $m_j,$ $j=1,2,3,$ and
let $f:M_1 \times M_2 \to M_3$ be a ${\cal C}^k$ map. The map $f$
is said to be differentiable (or smooth) with $k$-regularity in
$M_1$ if, for any charts
$\Big(U_1 \times U_2,
\big(x\equiv(x_1,\ldots,x_{m_1}),y\equiv(y_1,\ldots,y_{m_2})\big)\Big)\;
\mbox{in}\; M_1 \times M_2$
and $(U_3,z\equiv (z_1,  \ldots,z_{m_3}))\; \mbox{in}\; M_3,$
the local expression of $f,$ $f(x,y):x(U_1)\times y(U_2)
\to z(U_3),$ satisfies that $f(\cdot,y)$ is smooth in $x(U_1)$ for
any $y \in y(U_2),$ and all the partial derivatives of $f(x,y)$
with respect to variables in $x$ are ${\cal C}^k$ in $x(U_1)
\times y(U_2).$
\end{definition}

\begin{definition}[smooth deformation of the double of a circular domain]\label{def:smoothdef}
Let $v_0\in {\cal T}_n$ and $\epsilon >0$ small enough. Denote by
$V(\epsilon)$ the Euclidean ball of radius $\epsilon$ in ${\T}_n$
centered at $v_0.$ Since $V(\epsilon)$ is simply connected,
standard homotopy arguments in differential topology show the
existence of a family of diffeomorphisms $\{F_v:N(v_0) \to
N(v)\;:\; v \in V(\epsilon)\}$ such that $F_{v_0}=\mbox{Id},$
$F_v(\infty)=\infty,$ $J_v\circ F_v\circ J_{v_0}=F_v,$ for any $v
\in V(\epsilon),$ and  $F:V(\epsilon)\times
\overline{\Omega(v_0)}\to \c,$ $F(v,z):=F_v(z),$ is smooth.

By definition, we say that $\{F_v:N(v_0) \to N(v)\;:\; v \in
V(\epsilon)\}$ is a smooth deformation of $N(v_0).$ Moreover note
that, for $\epsilon$ small enough,  $\frac{\partial F}{\partial
z}\neq 0$ in $V(\epsilon)\times \overline{\Omega(v_0)}.$
\end{definition}

Let $W \subset \div_{k_1,k_2}$ be a submanifold, and let  $h: W
\to {\cal M}_n$ be a map preserving the fibers, that is to say,
$h_{v,D_1,D_2}:=h(v,D_1,D_2)\in{\cal M}(v)$ for any
$(v,D_1,D_2)\in {W}.$ We are going to define the notion of
differentiability with $k$-regularity of $h.$ Take ${\cal V}$ any
coordinate neighborhood in $\div_{k_1,k_2}$ meeting ${W}.$ Denote
by $V$ the associated neighborhood to ${\cal V}$ in ${\cal T}_n$
and call $v_0 \in V$ an interior point. Take a smooth deformation
of $N(v_0),$ $ \{F_v:N(v_0) \to N(v)\;:\; v \in V(\epsilon)\}.$ We
say that $h$ is {\em differentiable with $k$-regularity in ${\cal
V}\cap W$}  if the map
$$\hat{h}:({\cal V} \cap W) \times
N(v_0)\to \overline{\c},$$
$$\hat{h}((v,D_1,D_2),x)=h_{v,D_1,D_2}(F_v(x))$$
is smooth with
$k$-regularity in ${\cal V}\cap W$.  The map $h$ is said to be
differentiable with $k$-regularity on ${W}$  if it does in ${\cal
V}\cap W,$ for any coordinate neighborhood ${\cal V}$ meeting
${W}.$ It is easy to check that this definition does not depend on
choice of neither $v_0$ nor the smooth deformation of $N(v_0).$

Likewise, for a map $\omega:{W}\to {\cal H}_n$ preserving the
fibers, define $\hat{\omega}: {\cal V}({\epsilon})\cap W\to {\cal
H}(v_0)$ by
$$\hat{\omega}(v,D_1,D_2)=\big(
F_v^*(\omega_{v,D_1,D_2}) \big)^{(1,0)}$$
where the superscript
$(1,0)$ means the $(1,0)$ part of the 1-form (by definition
$(f\,dz+g\,d\overline{z})^{(1,0)} =f\,dz$). We say that $\omega$
{\em is differentiable  with
$k$-regularity in ${\V}(\epsilon)\cap W$} if for any local chart $(U,z)$ in $N(v_0),$ the
map $\hat{f}:({\V}({\epsilon})\cap W)\times U\to \overline{\c},$
given by $\hat{f}((v,D_1,D_2),z)=\hat{\omega}(v,D_1,D_2)(z)/dz$ is
smooth with $k$-regularity in ${\V}(\epsilon)\cap W.$ The global
concept of
differentiability with $k$-regularity in ${W}$ is defined in the obvious way.\\


The following 1-forms we will play an important role during the proof
of the Main Theorem.

Given $D=\prod_{j=1}^s w_j^{m_j} \in Div_k(\Omega(v)),$ we denote
by $\tau_D(v)$ the unique meromorphic 1-form on $N(v)$ satisfying:
\begin{itemize}
\item $\tau_D(v)$ has simple poles at $w_j$ and $J_v(w_j),$ $j=1,\ldots,s,$ and no other
poles,
\item $\mbox{Residue}_{w_j} \big(\tau_D(v)\big)= -\mbox{Residue}_{J_v(w_j)}\big(\tau_D(v)\big)=-m_j$ for any $j$
\item
$\int_{a_i(v)} \tau_D(v)=0,$ for any $i=1,\ldots,n.$
\end{itemize}

Likewise, take $D_1 =\prod_{j=1}^s w_{j,1}^{m_j},D_2=\prod_{h=1}^r
w_{h,2}^{n_h} \in Div_k(\Omega(v))$ and define
$\kappa_{D_1,D_2}(v)$ as the unique meromorphic 1-form on $N(v)$
satisfying
\begin{itemize}
\item $\kappa_{D_1,D_2}(v)$ has simple poles at $w_{j,1},\,w_{h,2}$ and
$J_v(w_{j,1}),\, J_v(w_{h,2}(v),$ $j=1,\ldots,s,$ $h=1,\ldots,r,$
and no other poles,
\item
$\mbox{Residue}_{w_{j,1}}\big(\kappa_{D_1,D_2}(v)\big)=\mbox{Residue}_{J_v(w_{j,1})} \big(\kappa_{D_1,D_2}(v)\big)=-m_j,$
for any $j$
\item
$ \mbox{Residue}_{w_{h,2}}\big(\kappa_{D_1,D_2}(v)\big)=\mbox{Residue}_{J_v(w_{h,2})}\big(\kappa_{D_1,D_2}(v)\big)=n_h$
for any $h$
\item  $\int_{a_i(v)} \kappa_{D_1,D_2}(v)=0,$ for any $i=1,\ldots,n.$
\end{itemize}

\begin{lemma}[\cite{fer-lop-sou}] \label{co:smooth}
The maps $\eta_j:{\cal T}_n\to {\cal H}_n,$ $v\mapsto\eta_j(v),$
$\tau:\div_{k}\to{\cal H}_n,$ $(v,D)\mapsto\tau_D(v),$ and
$\kappa:\div_{k,k}\to{\cal H}_n,$
$(v,D_1,D_2)\mapsto\kappa_{D_1,D_2}(v)$ are differentiable with
$1$-regularity.

As a consequence, the functions $\pi_{j,k}(v):= \int_{b_j(v)}
\eta_k(v),$ are differentiable on ${\cal T}_n.$
\end{lemma}

The following theorem will show that $\Mg_n$ and $\Gg_n$ are
analytic manifolds of dimension $3n+4.$ We first need the
following lemma, proved in \cite{fer-lop-sou}:

\begin{lemma}[\cite{fer-lop-sou}] \label{lem:porfin}
Given $v \in {{\cal T}_n},$ there exists a holomorphic 1-form
$\omega_0$ in $N(v)$ having $2n-2$ distinct zeroes, none of them
contained in $\partial \Omega(v),$ and satisfying $J_v^*
(\omega_0)=\overline{\omega_0}.$
\end{lemma}

\begin{theorem}[Main theorem] \label{th:structure}
The map
$$\sg_2:\Mg_n \to \big(\coc\big)^{n+1}\times]-1,1[, \qquad \sg_2(G,\mg)=(\mg,c),$$
where $c$ is the normal direction at the non-normalized end, is
injective and smooth. Hence, $\sg_2(\Mg_n)$ is open and so $\sg_2$ provides a global system of
analytic coordinates on $\Mg_n.$

Moreover, the action of the group of permutations of order $n+1,$
 $\mu:\Mg_n\times{\cal P}_n\to\Mg_n,$ is
discontinuous. Hence the orbit space, naturally identified to
$\Gg_n,$ has a unique analytic structure making
$\sg_1:\Mg_n\to\Gg_n,$ $\sg_1(G,\mg)=G,$ an analytic covering of
$(n+1)!$ sheets.
\end{theorem}

\begin{proof}

To see that $\sg_2$ is one to one, suppose $(G_i,\mg_i)\in\Mg_n,$
$i=1,2$ satisfy $\sg_1(G_1,\mg_1)=\sg_1(G_2,\mg_2)=(\mg,c).$ From
our normalizations, one end of both of them is asymptotic to
$\{x_3=0,x_2\geq 0\}.$ Since $G_1$ and $G_2$ are graphs over
$\{x_3=0\}\subset\coc$ it follows that, for both surfaces, the
other end is asymptotic to $\Pi\cap\{x_2\leq 0\},$ where $\Pi$ is
the plane determined by the normal direction $c.$ Therefore, $G_1$
and $G_2$ are contained in a common horizontal half space and
by Theorem \ref{th:uniqueness} we get $G_1=G_2.$\\

To finish the first part of the theorem,
it is enough to check that $\sg_2$ is smooth and then use the Domain Invariance Theorem. Here we have endowed $\Mg_n$
with the differentiable structure induced by ${\mathcal{E}}:\Mg_n\to\sb_n\times\coc\times\{1,-1\}$
(see Definition \ref{def:cale}).

Let $(S,\mg)\in\Mg_n$ and label $\chi=(x,q_0,\varepsilon_0)={\mathcal{E}}(S,\mg).$
Following the notation in the proof of the
proposition \ref{pro:cale}, call $X_\chi$ the associated maximal
immersion and label as $(g_\chi,\phi_3(\chi))$ its Weierstrass data.

\begin{quote}{ {\bf Claim} {\em  The
maps

$\begin{array}{ccc}
\sb_n\times\coc\times\{1,-1\}\to {\cal M}_n, & and &\sb_n\times\coc\times\{1,-1\} \to {\cal H}_n,\\
\chi \mapsto g_\chi & &\chi \mapsto \phi_3(\chi)
\end{array}$

are smooth with
$2$-regularity and $1$-regularity, respectively. Consequently, the
map $\sg_2$ is smooth.}}
\end{quote}

First we prove that the maps ${\cal S}_n \to {\cal M}_n,$ $ x
\mapsto g_x^0,$ and ${\cal S}_n \to {\cal H}_n,$ $x \mapsto
\phi_3^0(x),$ given by equations (\ref{eq:g}) and (\ref{eq:phi}),
are smooth with $2$-regularity and $1$-regularity, respectively.

Indeed, take $x_0=(v_0,D_0) \in {\cal S}_n.$ From Theorem
\ref{th:submanifold}, there exists an open ball $V(\epsilon)$ in
${\cal T}_n$  centered at $v_0$ of radius $\epsilon >0$ and a
local diffeomorphism ${V}(\epsilon) \to {\cal S}_n,$ $v \mapsto
(v,D(v)),$ where $D(v_0)=D_0.$ We label ${\cal V}(\epsilon)$ as
the image of $V(\epsilon)$ under this map. For simplicity, we
write $x(v):=(v,D(v)),$ $v \in { V}(\epsilon).$

Therefore,  the map ${V}(\epsilon) \to \div_{n+1},$ $v \to
(v,D(v)\cdot 0)$ is smooth, and since $\tau:\div_{n+1} \to {\cal
H}_n$ is also smooth with $1$-regularity (see Corollary
\ref{co:smooth}), the same holds for the map ${V}(\epsilon) \to
{\cal H}_n,$ $v \mapsto \tau_{v}:=\tau_{D(v)\cdot 0}(v).$

Take a smooth deformation of $N(v_0),$ $\{F_v:N(v_0)\to N(v)\;:\;v
\in V(\epsilon)\}.$  Let
$B(v_0)=\{a_1(v_0),\ldots,a_n(v_0),b_1(v_0),\ldots,b_n(v_0)\}$ be
the canonical homology basis on $\Omega(v_0)$ defined as in
Subsection \ref{sub:ident}. 
Observe that $N(v_0)-\cup_{j=1}^n (a_j(v_0) \cup b_j(v_0))$ is
simply connected, moreover, without loss of generality we can
suppose that this domain does not contains the points in $
D_0\cdot 0.$ For $v$ close enough to $v_0$ the curves
$a_j(v):=F_v(a_j(v_0)),$ $b_j(v_0):=F_v(b_j(v_0))$ are a canonical
basis of $N(v)$ and do not pass also through the points in $
D(v)\cdot 0,$ $j=1,\ldots,n.$

By Abel's theorem, and for $z \in N(v)-\cup_{j=1}^n (a_j(v) \cup
b_j(v))$ we have
$$g_{x(v)}^0(z)=\mbox{Exp} \Big( \int_1^z (\tau_{v}+ \sum_{j=1}^n m_j(v)
\eta_j(v))  \Big)$$
 In this expresion, the integration paths lie
in $\big(N(v)-\cup_{j=1}^n (a_j(v) \cup b_j(v))\big) \cup\{1\},$
and $m_j(v) \in \z$ are integer numbers determined by the
equation:
$$\widetilde{\varphi_v}(D(v)\cdot 0)-\widetilde{\varphi_v}(J_v(0) \cdot J_v(D(v)))=\sum_{j=1}^n
m_j(v) \pi^j(v),$$ where $\widetilde{\varphi_v}$ is the branch of
$\varphi_v$ on $N(v)-\cup_{j=1}^n (a_j(v) \cup b_j(v))$ vanishing
at $1.$

Since $m_j(v)$ depend continuously on $v,$ then $m_j(v)=m_j \in
\z$ and so, by Corollary \ref{co:smooth}, $g_{x(v)}^0$ depends
smoothly on $v$ with $2$-regularity.\\

We have to obtain the analogous result for the map ${\cal
V}(\epsilon) \to {\cal H}_n,$ $v \to \phi_3^0(x(v)).$ Take the
holomorphic 1-form $\omega_0$ on $N(v_0)$ given in Lemma
\ref{lem:porfin}, write $\nu(v_0):=\omega_0=\sum_{j=1}^n \lambda_j
\eta_j(v_0),$ where $\lambda_j \in \r,$ and define
$\nu(v):=\sum_{j=1}^n \lambda_j \eta_j(v).$ Since the map
$v\mapsto \nu(v)$ is smooth with $1$-regularity (see Corollary
\ref{co:smooth}) it suffices to prove that $v\mapsto
\frac{\phi_3^0(x(v))}{\nu(v)}$ is smooth with $2$-regularity.

By Hurwitz's Theorem and the implicit function theorem, $\nu(v)$
satisfies also the thesis in Lemma \ref{lem:porfin}, for $v\in
V(\epsilon),$ for $\epsilon>0$ small enough. Moreover, as
explained during the proof of Lemma \ref{lem:porfin}, the map
$V(\epsilon) \to \div_{2n-2},$ $v \mapsto (v,(\nu(v)))$ is at
least ${\cal C}^1,$ where as usually $(\nu(v))$ is the canonical
divisor associated to $\nu(v).$ Hence, writing $(\nu(v))=A_v \cdot
J_v(A_v),$ the map $V(\epsilon) \to \div_{n-1},$ $v \mapsto
(v,A_v),$ is also smooth, and therefore, the same holds for
${V}(\epsilon) \to \div_{n,n},$ $v \mapsto (v,D(v),\infty\cdot
A_v).$ We infer from Corollary \ref{co:smooth} that the  map
${V}(\epsilon)\to {\cal H}_n,$ $v \mapsto
\kappa_{v}:=\kappa_{\infty \cdot A_v,D(v)}(v),$ is smooth with
1-regularity. Reasoning as above, the map
$$f_{x(v)}(z)= \mbox{Exp} \Big( \int_1^z (\kappa_{v}+ \sum_{j=1}^n n_j
\eta_j(v))  \Big),$$ is a well defined meromorphic function on
$N(v),$ for suitable integer numbers $n_j$  not depending on $v$
and ${V}(\epsilon)\div_n \to {\cal M}_n,$ $v \mapsto f_{x(v)},$ is
smooth with $2$-regularity. The principal divisor associated to
$f_{x(v)}$ is given by $(f_{x(v)})=\frac{D(v) \cdot
J_v(D(v))}{\infty \cdot A_v \cdot J_v(\infty) \cdot J_v(A_v)}.$
Therefore, if we write $\nu(v)=h_v(z) \frac{dz}{z}$ on
$U(v)=\big(\Omega(v)-\{0,\infty\}\big)\cup\big(\Omega(v)^*-\{J_v(0),J_v(\infty)\}\big)\cup
a_0(v),$ we infer that
$\frac{\phi_3^0(x(v))}{\nu(v)}=\frac{1}{h_v(1)} f_{x(v)},$ and so
$v\mapsto \phi_3^0(x(v))$ is smooth with $1$-regularity.\\

It follows that the  map
$$X_x^0:=\mbox{Real}\int_1 \Phi(x)^0 \qquad
\Phi(x)^0:= \big(\frac{i}{2}(\frac{1}{g_x^0}- g_x^0),-\frac{1}{2}(\frac{1}{g_x^0}+ g_x^0), 1\big)\phi_3^0(x) $$
depends smoothly on $x=(v,D)$ with $2$-regularity and defines a
complete maximal surface in $\l^3/_{\langle V_x\rangle},$ where
$$V_x=\mbox {Real}\big[2\pi i \mbox{ Res}\,_0 (\Phi(x)^0) \big] $$
Since $\phi_3^0(x)$ is holomorphic at $0,$ the vector
$V_x=(w_x,0)\in \c\times \r$ is horizontal.
Moreover, up to
replacing $\phi_3^0(x)$ by $-\phi_3^0(x)$ we can suppose that the
third coordinate of the flux around the curve $a_0(v)$ is
positive.

It is straightforward to see that $\theta_\chi= \frac{\bar
{w}_x}{|w_x|}$ and $r_\chi =\frac{\varepsilon_0}{|w_x|}$ depend
smoothly on $\chi,$  and so, it follows that $g_\chi=\theta_\chi
g_x^0$ and $\phi_3(\chi)=r_\chi\phi_3^0(x) $ depend smoothly on
$\chi$ with 2 and 1-regularity respectively.

To conclude the proof of the claim, observe that
$$X_\chi=q_0+   \mbox{Real}\int_1\big(\frac{i}{2}(\frac{1}{g_\chi}-
g_\chi),-\frac{1}{2}(\frac{1}{g_\chi}+ g_\chi), 1\big)\phi_3^0(\chi)  $$

depends smoothly on $\chi=(x=(v,D),q_0,\varepsilon_0)$ with 2-regularity. Therefore
$q_j(\chi)=X_\chi(a_j(v))$ and $c(\chi)=g_\chi(\infty)$ are smooth  functions of
$\chi,$ and  the same holds for $\sg_2.$ This proves the
claim.

\vspace{0.3cm}

By the injectivity of $\sg_2$ and the domain invariance theorem,
${\sg}_2(\Mg_n)$ is an  open domain in $\big(\coc\big)^{n+1} \times  ]-1,+1[.$
We can
then endow $\Mg_n$ with the unique analytic structure making
${\sg}_2: {\Mg}_n \to
{\sg}_2({\Mg}_n)$ an analytic diffeomorphism.\\

To conclude,  it remains to check that the action $\mu$ is
discontinuous. Indeed, let $\tau:\Mg_n \to \Mg_n$ denote the
diffeomorphism given by $\tau(S,\mg)=(S,\tau(\mg)),$ $\tau \in
{\cal P}_{n+1}.$ Let $(S_0,\mg_0) \in \Mg_n$ and write
$\mg_0=(q_0,q_1,\ldots,q_n) \in \big(\coc\big)^{n+1} .$ Take a
neighborhood $U_j$ of $q_j$ in $\l^3/_{\langle(1,0,0)\rangle},$
$j=0,1,\ldots,n,$ such that $U_i \cap U_j=\emptyset,$ $i \neq j,$
and call ${\cal U}=\prod_{j=0}^n U_j.$  Then, it is clear that
$\tau(\sg_2^{-1}({\cal U}\times \r)) \cap \sg_2^{-1}({\cal
U}\times \r)=\emptyset,$ for any $\tau \in {\cal
P}_{n+1}-\{\mbox{Id}\},$ which proves the discontinuity of $\mu$
and concludes the proof.

\end{proof}
To finish, we prove that the underlying topology in $\Gg_n$ corresponds to the
uniform convergence of graphs over compacts subsets of $\{x_3=0\}.$
\begin{theorem}
Let $\{G_k\}_{k \in \n}$ be a sequence in $\Gg_n,$ and $G_0 \in \Gg_n.$

Then $\{G_k)\}_{k \in \n} \to G_0$ in the topology of $\Gg_n$ if and only if
$\{G_k\}_{k \in \n}$ converges to $G_0$ uniformly on compact subsets of $\{x_3=0\}.$
\end{theorem}
\begin{proof}
Suppose $\{G_k\}_{k \in \n} \to G_0 \in \Gg_n$ in the topology of $\Gg_n,$ and
choose marks in such a way that $\{(G_k,\mg_k)\}_{k \in \n}$  converges  to  $(G_0,\mg_0) $ in $\Mg_n.$

Write ${\cal E}((\Gg_k,\mg_k))=(x_k,q_0(k),\epsilon_k),$ $X_k=X_{(x_k,q_0(k),\epsilon_k)}$ and $x_k=(v_k,D_k) \in {\cal S}_n,$
$k \in \n \cup \{0\}.$ Observe that, without loss of generality, $\epsilon_k=\epsilon_0,$ for all $k \in \n.$

Since $\{(x_k,q_0(k),\epsilon_k)\} \to (x_0,q_0(0),\epsilon_0)$ and $X_{(x,q_0,\epsilon)}$ depends smoothly on $(x,q_0,\epsilon)$ with 2-regularity (see the proof of Theorem \ref{th:structure}),  it is not hard to check that $\{X_k\}_{k \in \N}$ diverges uniformly on $k,$ that is to say, for any compact $W$ in $\{x_3=0\}$ there is $r>0$ such that $|z+1/z|>r$ implies $X_k(z) \notin W \times \r,$ for all $k.$

Let $W$ be any compact domain in the cylinder $\{x_3=0\} \subset \coc $ containing the singularities in $\mg_0$ as interior points, and let $W_k$ denote the compact set $X_k^{-1}(W \times \r) \subset \overline{\Omega(v_k)}-\{0,\infty\},$ $k \in \n \cup \{0\}.$

As the domains $W_k$ are uniformly contained  in a compact region of $\c-\{0\},$ then $\{W_k\}_{k \in \n} \to W_0$ in the Hausdorff distance and $X_k$ converges uniformly on $W_0$ to $X_0.$ In the last statement we have used that $X_k$ can be reflected analytically about the circles in $\partial {\Omega}(v_k),$ and so all the immersions $X_k,$ $k$ large enough, are well defined in a universal neighborhood of $W_0$ in $\c.$ It is then obvious that the function $u_k: \r^2 \to \r$ defining the graph $G_k$ converges uniformly over $W$ to the function $u_0: \r^2 \to \r$ defining $G_0$ (furthermore, $\{v_k\}_{k \in \n} \to v_0$ implies that
$\{\mg_k\}_{k \in \n} \to \mg_0$). Since $W$ can be as larger as we want, $\{u_k\}_{k \in \n} \to u_0$ uniformly on compact subsets of $\r^2.$

Assume now that the functions $u_k:\{x_3=0\} \to \r$ defining $G_k$ converge, as $k \to \infty,$ to the function $u_0:\{x_3=0\} \to \r$ defining  $G_0$ uniformly on compact subsets of $\{x_3=0\}.$

Let us show that  singular points of $G_0$ are  limits of sequences of singular points of graphs $G_k,$ $k \in \n.$
Indeed, let $p_0=(y_0,u_0(y_0)) \in G_0$ be a singular point, and without loss of generality, suppose that $p_0$ is a downward pointing conelike singularity. By Lemma \ref{lem:sing}, there exists $\epsilon>0$ small enough such that $u_0^{-1}(\{x_3 \leq u_0(y_0)+\epsilon\})$ contains a compact component $C_0(\epsilon)$ with regular boundary and containing $y_0$ as the unique (interior) singular point.
Since $\{u_k\}_{k \in \n} \to u_0$ uniformly on compact subsets, $u_k^{-1}(\{x_3 \leq u_0(y_0)+\epsilon\})$
must contain a compact component $C_k(\epsilon)$ containing $y_0$ as well, $k$ large enough. Furthermore, $\{C_k(\epsilon)\} \to C_0(\epsilon)$ in the Hausdorff sense, and by the maximum principle $C_k(\epsilon)$ must contain at least an interior singular point $y_k$ of $u_k,$ $k$ large enough. Since $C_0(\epsilon)$ converges to $\{y_0\}$ as $\epsilon \to 0,$ we deduce that $\{p_k:=(y_k,u_k(y_k))\}_{k \to \infty} \to p_0.$

As a consequence, there exist marked graphs $(G_k,\mg_k) \in \Mg_n,$ $k \in \n \cup \{0\},$ such that $\{\mg_k\}_{k \in \n} \to \mg_0.$

Call $c_k$ the stereographic projection of the normal vector at $E_2$ of $G_k,$ $k
\in \n \cup \{0\},$ and let us see that $\{c_k\} \to c_0.$ Indeed, take $C<0$ small
enough such that $\{x_3=0, x_2 \leq C\}$ does not contain any singular point of
$u_k,$ $k \in \n \cup \{0\},$  let $\gamma$ be a circle in $\{x_3=0, x_2 \leq C\}$
and let $A$ denote a closed tubular neighborhood of $\gamma$ in $\{x_3=0\}$ not
containing any singular point of $u_0.$ It is well known that  $u_k-u_0$ is solution
of a uniformly elliptic linear  equation $L_k(u_k-u_0)=0$ over $A,$ $k$ large
enough. Moreover, the fact that the functions $\frac{1}{1-|\nabla u_k|},$ $k \in
\n,$ are uniformly bounded on $A$ (see \cite{Bar-Sim}) guarantee that the
coefficients of operators $L_k,$ $k \in \n,$ are uniformly bounded too. Therefore,
since $\{u_k\}_{k \in \n} \to u_0$ uniformly on $A,$   the classical Schauder
estimates (\cite{gilbarg} p. 93) imply that $\{u_k\}_{k \in \n} \to u_0$ in the
$C^2$ norm on $A.$ In particular, $$\{n_k:= \int_\gamma \nu_k(s_k) ds_k  \}_{k \in
\n} \to n_0:=\int_\gamma \nu_0(s_0) ds_0,$$ where $\nu_k$ and $s_k$ are the conormal
vector and the arc-length parameter along $\gamma$  in $G_k,$ respectively, for any
$k \in \n \cup \{0\}.$ Since the normal vector at $E_2$ of $G_k$ lies in $\{x_1=0\}$
and is orthogonal to $n_k,$ $k \in \n\cup\{0\},$  we infer that $\{c_k\} \to c_0.$

Since $\sg_2:\Mg_n \to \sg_2(\Mg_n) \subset \r^{3n+4}$ is an homeomorphism,  $\{(G_k,\mg_k)\}_{k \in \n} \to (G_0,\mg_0)$ in the manifold $\Mg_n,$ and so,
$\{G_k\}_{k \in \n} \to G_0$ in the manifold $\Gg_n.$ This concludes the proof.
\end{proof}

{\bf ISABEL FERNANDEZ, FRANCISCO J. LOPEZ, \newline Departamento
de Geometría y Topología \newline Facultad de Ciencias,
Universidad de Granada \newline 18071 - GRANADA (SPAIN) \newline
e-mail:(first author) isafer@ugr.es, (second author)
fjlopez@ugr.es

\vspace{0.2cm}

RABAH SOUAM, \newline Institut de Mathématiques de Jussieu-CNRS
UMR 7586\newline Université Paris 7\newline Case 7012\newline
2,place Jussieu\newline 75251 Paris Cedex 05, France
\newline
e-mail: souam@math.jussieu.fr}


\begin{thebibliography}{99}
\bibitem{Bar-Sim} R. Bartnik and L. Simon: {\em  Spacelike hypersurfaceswith prescribed boundary values and mean curvature.} Comm. Math. Phys.,Vol. 87(1982/83), 131-152.

\bibitem{calabi} E. Calabi.: {\em Examples of the Bernstein problem for
some nonlinear equations.} Proc. Symp. Pure Math., Vol. 15,
(1970), 223-230.

\bibitem{cheng-yau} S. Y. Cheng and S. T. Yau.: {\em  Maximal space-like
hypersurfaces in the Lorentz-Minkowski spaces.}  Ann. of Math.
(2), Vol. 104 (1976), 407-419.

\bibitem{ecker} K. Ecker.: {\em  Area maximizing hypersurfaces in
Minkowski space having an isolated singularity .}  Manuscripta
Math., Vol. 56 (1986), 375-397.

\bibitem{Es} F. J. M. Estudillo and A. Romero.: {\em Generalized maximal
surfaces in the Lorentz-Minkowski space $\l^3$} Math. Proc. Camb.
Phil. Soc. 111, (1992), 515-524.

\bibitem{farkas} H. M. Farkas, I. Kra.: {\em Riemann surfaces.} Graduate
Texts in Math., {\bf 72}, Springer Verlag, Berlin, 1980.

\bibitem{fer-lop} I. Fern\'{a}ndez and  F. J. L\'{o}pez.:
{\em Periodic Maximal surfaces in the Lorentz-Minkowski space $\l^3$.} Preprint.

\bibitem{fer-lop-sou} I. Fern\'{a}ndez, F. J. L\'{o}pez and R. Souam.:
{\em The space of complete embedded maximal surfaces with isolated
singularities in the 3-dimensional Lorentz-Minkowski space
$\l^3.$} Preprint. Link: http://arxiv.org/PS$\underline{\,\;}$cache/math/pdf/0311/0311330.pdf

\bibitem{gilbarg} D. Gilbarg and N. S. Trudinger.: {\em Elliptic partial
differential equations of second order.} Springer-Verlag, (1977).


\bibitem{klyachin} A. A. Klyachin.: {\em Description of a set of entire solutions
with singularities of the equation of maximal surfaces.}
(Russian) Mat. Sb. 194 (2003), no. 7, 83--104; translation in Sb. Math. 194 (2003), no. 7-8, 1035--1054.

\bibitem{tubes} V. A. Klyachin and V.M. Miklyukov.: {\em Geometric structures of tubes and bands of zero mean curvature in Minkowski space.} Annales Academia Scientiarum Fennicae Mathematica, 28 (2003) 239-270

\bibitem{kobayashi} O. Kobayashi.: {\em Maximal surfaces with conelike
singularities.} J. Math. Soc. Japan 36 (1984), no. 4, 609--617

\bibitem{lopez-lopez-souam} F. J. L\'{o}pez, R. L\'{o}pez and R. Souam.:
{\em Maximal surfaces of Riemann type in Lorentz-Minkowski space
$\l^3$.} Michigan J. of Math., Vol. 47 (2000),  469-497.

\bibitem{um-ya} M. Umehara and K. Yamada: {\em Maximal surfaces with singularities in Minkowski space.} Preprint.

\bibitem{wolf} J. Wolf: {\em Spaces of Constant curvature.} McGraw-Hill, New York (1967).

\bibitem{oneill} B. O'Neill: {\em Semmi-riemannian geometry.} Academic Press (1983).


\end{thebibliography}
\end{document}